\documentclass [10pt,reqno]{amsart}
\usepackage{ucs}
\usepackage[latin1]{inputenc} 
\usepackage{amsmath, amssymb, amscd}
\usepackage{pb-diagram}
\usepackage[latin1]{inputenc}
\usepackage{amsmath}
\usepackage{amssymb}
\usepackage[initials]{amsrefs}

\usepackage{color}

\usepackage{graphicx}	
\usepackage{picins}

\newtheorem{theorem}{Theorem}[section]
\newtheorem{definition}[theorem]{Definition}
\newtheorem{lemma}[theorem]{Lemma}
\newtheorem{prop}[theorem]{Proposition}
\newtheorem{corollary}[theorem]{Corollary}

\newtheorem{remark}[theorem]{Remark}
\newtheorem{corollary of proof}[theorem]{Corollary of Proof}

\newcommand{\abs}[1]{\lvert#1\rvert}

\newcommand{\modulo}[1]{\ ( \text{mod} \ #1 \, )}

\renewcommand{\epsilon}{\varepsilon}

\DeclareMathAlphabet{\mathpzc}{OT1}{pzc}{m}{it}
\usepackage{mathrsfs}

\newcommand{\CS}{\text{\em CS}}

\newcommand{\Z}{\mathbb{Z}}
\newcommand{\C}{\mathbb{C}}

\renewcommand{\H}{\mathbb{H}}
\newcommand{\R}{\mathbb{R}}

\renewcommand{\qed}{$\hfill \square$ \bigskip \\}
\renewcommand{\phi}{\varphi}

\newcommand{\Ad}{\text{Ad}}

\newcommand{\Hom}{\text{Hom}}

\newcommand{\id}{\text{id}}

\newcommand{\Map}{\text{Map}}
\newcommand{\sign}{\text{sign}}

\newcommand{\su}{\mathfrak{su}}

\renewcommand{\bar}{\overline}

\newcommand{\tr}{\text{tr}}

\newcommand{\zetadiez}{\zeta^{\#}}

\renewcommand{\i}{{\bf i}}
\renewcommand{\j}{{\bf j}}
\renewcommand{\k}{{\bf k}}
\newcommand{\x}{{\bf x}}
\newcommand{\y}{{\bf y}}
\newcommand{\z}{{\bf z}}

\renewcommand{\S}{ S(\operatorname{Im}(\H))}

\begin{document}
\thispagestyle{empty}
\title{Representation spaces of pretzel knots} 
\author[Raphael Zentner]{Raphael Zentner} 
\address {Mathematisches Institut \\
Westf\"alische Wilhelms-Universit\"at M\"unster \\
Einsteinstra\ss e 62 \\
48149 M\"unster
\\
Germany}
\email{raphael.zentner@math.uni-muenster.de}
\maketitle
\begin{abstract}
We study the representation spaces $R(K;\bf{i})$ as appearing in Kron\-heimer and Mrowka's  instanton knot Floer homologies, for a class of pretzel knots. In particular, for pretzel knots $P(p,q,r)$ with $p, q, r$ pairwise coprime, these appear to be non-degenerate and comprise representations in $SU(2)$ that are not binary dihedral. 
\end{abstract}

\section{Introduction}
Let $K$ be a knot in the 3-sphere, and $y_{0}$ a point in its complement. Let furthermore $m$ be a meridian of the knot. In the construction of framed instanton knot Floer homology \cite{KM_knots} there appear at the chain group level representation spaces 
\begin{equation*}
	R(K;\i)= \{ \rho \in \Hom (\pi_{1}(S^{3} \setminus K;y_{0}),SU(2))  \ | \ \rho(m) \sim \i \} 
\end{equation*}
of knots with the meridian $m$ (or links with each of their meridians) mapped to traceless matrices, or, equivalently, to elements that are conjugate to $\i$ when $SU(2)$ is viewed as the group of unit quaternions. In the construction of the slightly newer reduced singular knot Floer homology \cite{KM_sutures, KM_ss} there appears a very related representation space. 

Our intention is to study these representation spaces for a class of pretzel knots. In Section 3, we describe the conjugacy classes of representations of $P(p,q,r)$ by triangles on the 2-sphere, with two vertices fixed in order to fix the conjugacy class, and with the length of each of the edges taken from a finite set. The binary dihedral representations appear as degenerate triangles, where all three vertices lie on a great circle. More generally, we show that for knots $K=P(p_{1}, \dots, p_{n})$ the conjugacy classes of representations are described by (ordered) n-gons on the 2-sphere, again with two vertices fixed for fixation of the conjugacy class, and with the lenghts of the edges taken from a discrete set. As a consequence, each conjugacy class of representation of a 3-strand pretzel knot $P(p,q,r)$ is isolated in $\mathscr{R}(K;\i) := R(K;\i)/SU(2)$, and comes in general with an $n-3$ dimensional family for the case of $n$ strands. We emphasise this result in Section 5, where we compute the Zariski tangent spaces of all representations of $P(p,q,r)$ provided some arithmetic conditions on $p,q$ and $r$ hold. The required arithmetic conditions are studied in Section 4. In Section 6 we summarise the results for the pretzel knot $P(p,q,r)$, we consider examples, and we reobtain very easily the known result that a pretzel knot or link has bridge number 3 if $\abs{p}, \abs{q}, \abs{r} > 1$ under the assumption of pairwise coprimeness of these numbers. 
\newpage 

\subsection{Relation to Khovanov homology} 
There is an isomorphism of abelian groups
\begin{equation}\label{Kh=H(R)}
	\text{\em Kh}(K) \cong H_*(R(K;\i);\Z) \ ,
\end{equation}
where $\text{\em Kh}(K)$ denotes the Khovanov homology \cite{Kh} of $K$, for certain knots. For torus knots of type $(2,p)$ this was observed by Kronheimer and Mrowka \cite[Observation 1.1]{KM_knots}. For an arbitrary 2-bridge knot or 2-component link this was proved by Lewallen \cite{Lewallen}(in the current version by use of an unpublished result of Shumakovitch \cite{Shumakovitch}). More precisely, he shows that Khovanov homology of a one or two component alternating link is isomorphic to the integer homology of $R_{bd}(K;\i)$, where $R_{bd}(K;\i) \subseteq R(K;\i)$ is the subspace of binary dihedral representations. 
Our explicit description in Proposition \ref{summary for representation space} allows us to draw the following conclusion:
\begin{prop}\label{inequality}
Let $K$ be the alternating pretzel knot $P(p,q,r)$ for $p, q, r$ pairwise coprime, and such that $\abs{p},\abs{q},\abs{r} > 1$. Then 
\begin{equation}\label{Kh not equal H(R)}
	\text{\em Kh}(K) \ncong H_*(R(K;\i);\Z) \ ,
\end{equation}
i.e. these two abelian groups are not isomorphic. 
\end{prop}
{\em Proof:} In fact we have $R(K;\i) \cong R_{bd}(K;\i) \coprod \, (\coprod_{I'} \mathbb{RP}^3)$, where $I'$ parametrises the non-empty set of conjugacy classes of non-binary dihedral representations. 
\qed

\subsection{Relation to Kronheimer and Mrowka's instanton knot Floer homology}
In \cite{KM_knots} Kronheimer and Mrowka construct an abelian group called the framed instanton Floer homology $\text{\em FI}_*(K)$ of knots $K$ in the 3--sphere, and then later introduced an equivalent theory $I^\#(K)$ in \cite{KM_ss} that is there called `unreduced singular knot Floer homology' of $K$, and a reduced version $I^\natural(K)$. These are Morse homologies of a Chern-Simons functional $\text{\em CS}$ defined on a space of connections on an open 3-manifold obtained $Y$ from $K$. For the unreduced theory the 3-manifold $Y$ is the knot complement of $K$ with in addition an unlinked Hopf link and an arc joining the two components of the Hopf link removed. For the reduced version the 3-manifold $Y$ is the knot complement with a meridian and an arc going from the knot to the meridian removed. In the unreduced theory the critical space is just 
$
	R(K;\i) 
$
in our notation, and it is a related space $\overline{R}(K;\i)$ for the reduced theory, where each conjugacy class of an irreducible representation appears with a 1-sphere $S^1$, and where the reducible representation appears as a point. To be more precise, the requirement that a meridian is mapped to an element $\rho(m)$ that is conjugated to $\i$ is replaced with the requirement that $\rho(m) = \i$. The stabiliser $S^1$ acts freely on irreducible representations and trivially on abelian ones. 

In both cases, the critical space is degenerate, and so the Chern-Simons functional $\CS$ has to be perturbed. But if the critical space is non-degenerate in the Morse-Bott sense one expects the following to happen: The filtration coming from the Floer grading induces a spectral sequence starting at the homology of the critical manifold of $\CS$ and that converges to the instanton Floer homology of a generically perturbed functional. This attempt is chosen for instance by Fukaya \cite{Fukaya}, where the corresponding classical result in Morse homology is proved and then modified so as to be applicable to the instanton Floer homology of the connected sum of two manifolds, where the critical space is necessarily degenerate. The author has learned about this reference from Daniel Ruberman. 

Presumably our results of Section \ref{non deg} imply that for the pretzel knots $P(p,q,r)$ with $p,q,r$ coprime  the critical space 
\begin{equation*}
	\bar{R}(P(p,q,r);\i) 
\end{equation*}
is non-degenerate in the sense of \cite{KM_knots} in the normal directions, and so non-degenerate in the Morse-Bott sense for the setting of \cite{KM_knots}. If this is the case then our results of Proposition \ref{summary for representation space} indicate that this spectral sequence  has non-trivial differentials at least when all of $p,q,r$ have the same sign (which corresponds to an alternating pretzel knot). In fact, for alternating knots $K$ the rank of ${I}^\natural(K)$ is known \cite[Corollary 1.6]{KM_ss} to be equal to $\abs{\Delta_K(-1)}$, the absolute value of the knot determinant, whereas the rank of 
\[H_*(\bar{R}(P(p,q,r);\i);\Z) \] 
is strictly bigger if $\abs{p},\abs{q}$ and $\abs{r}$ are all strictly bigger than $1$. Therefore, a spectral sequence as made allusion to would have to have non-trivial differentials. 

For a 2-bridge knot, the ranks of the groups $I^\natural(K)$ and $H_*(\bar{R}(K;\i))$ are equal by Klassen's result \cite{Klassen}, determining the number of irreducible binary dihedral representations (and there are only such for a 2-bridge knot) to be equal to $(\abs{\Delta_K(-1)} - 1)/2$,  and the above cited result \cite[Corollary 1.6]{KM_ss}. This suggests that after generic perturbation of the Chern-Simons functional each $S^1$ coming from an irreducible representation gives rise to two critical points, the reducible representation to one critical point, and that the differentials are all trivial. For the alternating knots $P(p,q,r)$ this cannot happen by the above given rank argument.  
From the author's naive point of view this is surprising in the sense that it just looks like the non binary dihedral representations `do not contribute' to $I^\natural(P(p,q,r)$ for $p,q,r$ of the same sign, at least regarding the rank.

Thinking further, it would be interesting to compute the differentials yielding $I^\natural(K)$ or $I^\#(K)$  explicitly and to study Question 1.2 of \cite{KM_knots} explicitly on the class of pretzel knots considered here. 

\subsection{Parallels with representation spaces of Brieskorn homology spheres} \label{parallels}
There appear to be parallels between the representation spaces $R(K;\i)$ for $K=P(p_{1}, \dots, p_{n})$ and the representation spaces $R(Y) = \Hom(\pi_{1}(Y);SU(2))$ for $Y = \Sigma(a_1, \dots, a_n)$ a Seifert fibred homology sphere \cite{FS}. In both cases the representation space is non-degenerate for $n=3$ and degenerate for $n \geq 4$, with a similar growth in the dimensions of the Zariski tangent spaces. However, the analogies between the two cases also have limitations: In the case of the Brieskorn homology spheres the Floer gradings of the critical points all have the same parity, so there are no non-zero differentials in the instanton Floer chain complex, and the Floer homology is just isomorphic to the chain complex. As indicated above, a similar statement doesn't seem to be true for $I^\natural(P(p,q,r))$ for $p,q$ and $r$ pairwise coprime and of the same sign.  

The author was informed by Saveliev that the parallels on the critical space level can be explained by results of Collin and Saveliev, and Saveliev. For knots the representation spaces as considered here correspond to $\Z/2$ equivariant representation spaces of the manifold obtained as double branched cover of the knot complement, branched along the knot \cite[Proposition 3.3]{CollinSaveliev}. In the case one gets a Brieskorn homology sphere (an integer homology sphere), all representations turn out to be equivariant \cite{Saveliev_Brieskorn}, explaining the analogy we have encountered, and there is a study of the equivariant representations in the case one gets a more general Seifert fibred integer homology sphere in \cite{Saveliev_Seifertfibred}. Saveliev's result is applicable to more general Montesinos knots than just pretzel knots of which the branched cover is a Brieskorn homology sphere. The double branched cover of a knot in $S^3$ is an integer homology sphere if and only if the determinannt $\abs{\Delta_K(-1)}$ is one. 
Our result extends Saveliev's in the sense that it also yields an analogy in the case the double branched cover isn't an integer homology sphere, like for most pretzel knots. 
\\


\noindent In a former version of this article the consideration was restricted to alternating knots $P(p,q,r)$ with all of $p,q,r$ odd. When preparing a seminar talk on the topic the author realised that this restriction was unnecessary. 

\section*{Acknowledgements} The author wishes to thank Kim Fr\o yshov, Andrew Lobb, Stefan Friedl and an anonymous referee for their comments on preliminary versions of this manuscript. He would also like to thank Nikolai Saveliev and Daniel Ruberman for comments on the first arxiv-version of this article. Furthermore, he is grateful to Wolfgang L\"uck for financial support through his Leibniz prize. 

\section{Presentations of pretzel knot groups}

Let us consider the elementary `$p$-tangle' as in the figure. This shall mean  that we have a braid with $\abs{p}$ crossings, that is a `left-hand screw' if $p$ is positive, as in our figure, and that is a `right-hand screw' for $p$ negative. Let $s$ and $t$ be meridians at the top as indicated (with the basepoint in front of the eye of the observer and straight lines going directly to the starting point and end point of the indicated flash), and $u$ and $v$ at the bottom. 

In the situation of the figure, that is, for a left-handed $p$-tangle (with $p$ positive) we get 
\parpic[r]{\includegraphics[width=0.3\textwidth]{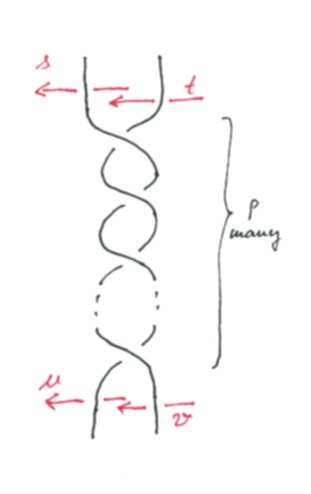}}
\noindent
\begin{equation*}
	\begin{split}
	u &  = (ts)^{-k}\ s^{-1} t s \ (ts)^{k}  , \\
	v &  = (ts)^{-k}\ s \ (ts)^{k}
	\end{split}
\end{equation*}
for $p=2k+1$ odd, and 
\begin{equation*}
	\begin{split}
	u &  = (ts)^{-k}\ s \ (ts)^{k} , \\
    v &  = (ts)^{-(k-1)}\ s^{-1} t s  \ (ts)^{k-1}  
	\end{split}
\end{equation*}
for $p=2k$ even.
On the other hand, if $p$ is negative, then we get
\begin{equation*}
	\begin{split}
	u &  = (ts)^{k}\ t \ (ts)^{-k}  , \\
	v &  = (ts)^{k}\ t s t^{-1} \ (ts)^{-k}
	\end{split}
\end{equation*}
for $p=-(2k+1)$ odd, and 
\begin{equation*}
	\begin{split}
	u &  = (ts)^{k-1}\ t s t^{-1} \ (ts)^{-(k-1)} , \\
    v &  = (ts)^{k}\ t \ (ts)^{-k}
	\end{split}
\end{equation*}
for $p=-(2k)$ even.
\\

With these formulae at hand we get the complement of any pretzel knot rather quickly. Indeed, the figures below show a diagram of a general pretzel knot or link $P(p_{1}, \dots, p_{n})$, where the boxes are to be filled in by the elementary tangles as in the figure above, and as a concrete example that of the knot $P(-3,5,7)$. 
\begin{figure}[ht]
\includegraphics[width=0.6\textwidth]{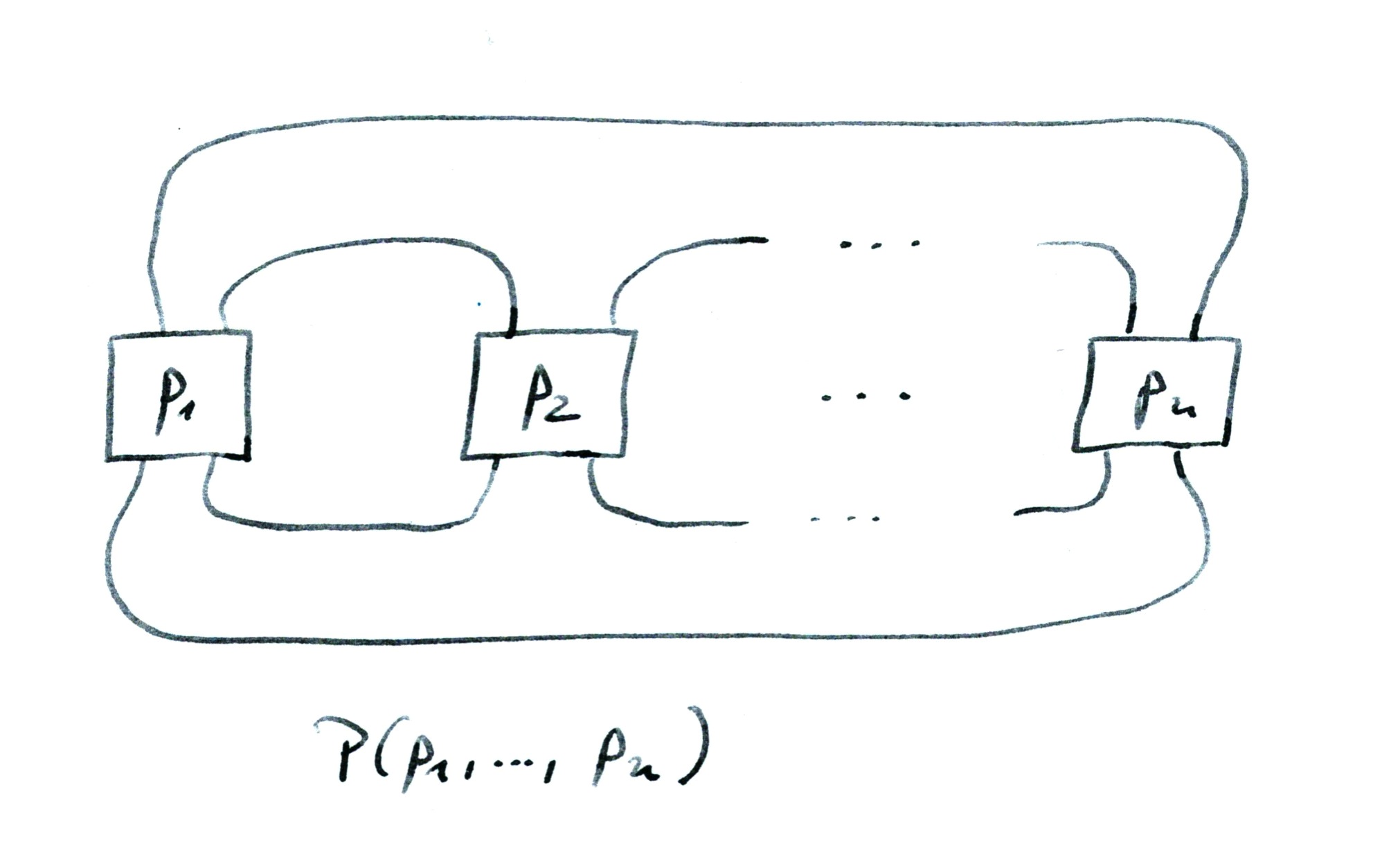}
\end{figure}


In this diagram it becomes obvious that it can be visualised as a knot or link with $n$ bridges or a  $2n$--plat. So there is an embedded 2-sphere (the `horizontal' one) in the 3-sphere cutting the $P(p_{1}, \dots, p_{n})$ pretzel knot/link in $2n$ points such that the resulting balls each contain $2n$ unknotted arcs with boundaries on the boundary 2-sphere. We therefore see that the knot or link complement has a decomposition into two pieces which each are (or deformation retract onto) two handlebodies of genus $n$ with common intersection a two-sphere punctured in $2n$ discs, with these discs centered at the intersection points of the knot or link with this 2-sphere. 
\parpic[r]{\includegraphics[width=0.45\textwidth]{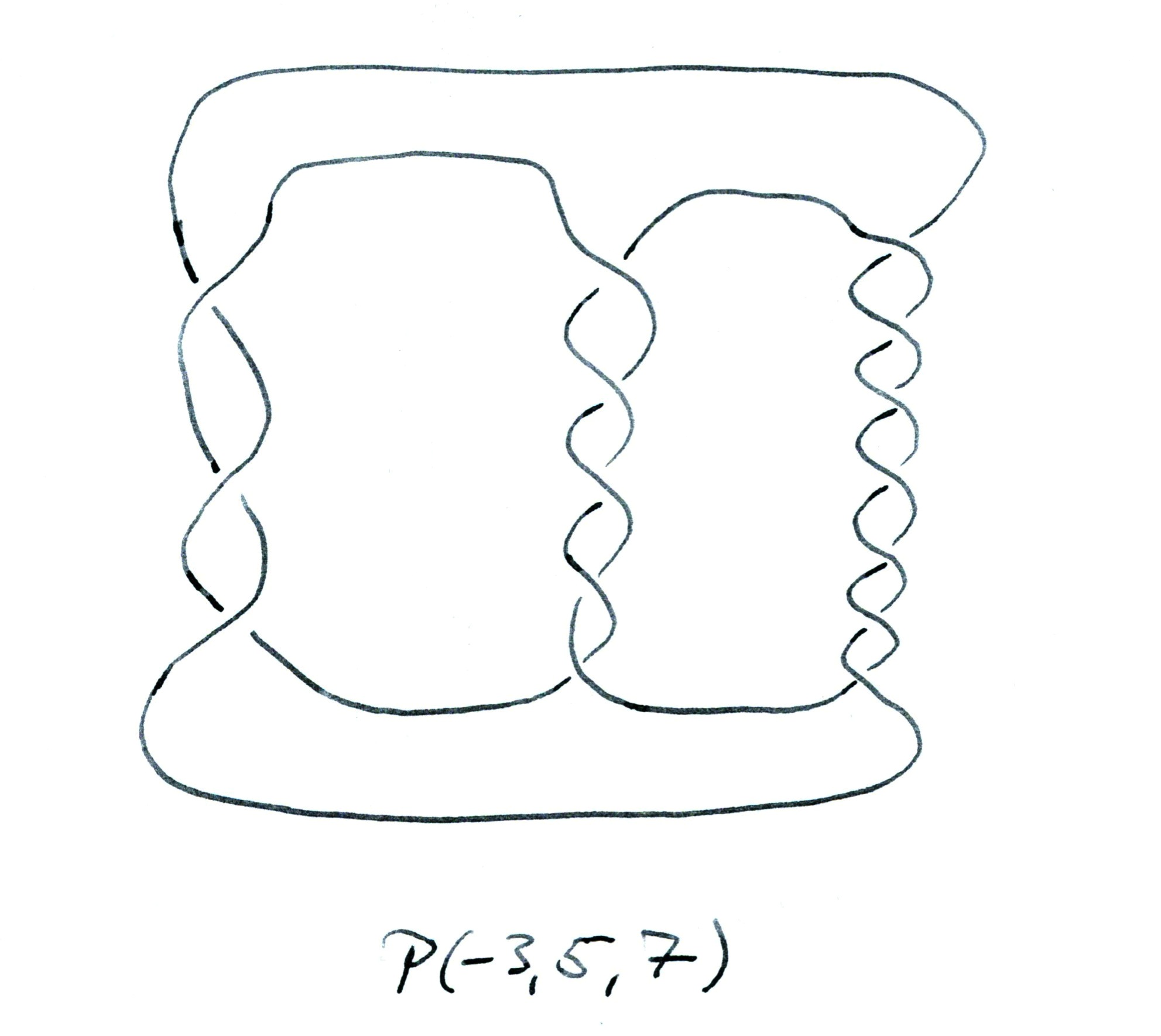}}
The fundamental group of each handlebody is a free group on $n$ generators given by the meridians of the knot/link in the corresponding ball. 

Now let $s_{1}, \dots, s_{n}$ be these meridians in the upper and $u_{1}, \dots, u_{n}$ be the merdians in the lower handlebody. We orient these meridians so that at each `elementary p-tangle' at the position $i$ we have
$s_{i}$ corresponding to $s$ in the picture above, and $s_{i+1}^{-1}$ corresponding to $t$ in the picture above. 

The Seifert-van-Kampen theorem now gives a presentation of the knot/link complement in a straight-forward manner. 

\begin{prop}\label{presentation knot group}
The fundamental group $G(K):= \pi_1(S^{3} - P(p_1,\dots,p_n),y_0)$ of the complement of the pretzel knot or link $P(p_{1}, \dots, p_{n})$  is given by
\begin{equation*}
\begin{split}
	\langle s_{1}, \dots, s_{n} | \, 
	v_1 u_2 = 1 \ ,
	v_2 u_3 = 1 \ , 
	\dots \ ,
	v_{n-1} u_n = 1 \ ,
	v_n u_1 = 1
	\rangle  ,
\end{split}
\end{equation*}
where $v_k = u_{k+1}^{-1}$ for $k=1, \dots, n-1$, where $v_n = u_1^{-1}$, and where $u_k$ and $v_k$ are the meridians at the bottom of the $k^{\text th}$ elementary tangle, as indicated in the figure, and as given in the formulae above in terms of $s_k$ and $t_k = s_{k+1}^{-1}$.
Each generator is a meridian, and any relation is a consequence of all others, so (any) one relation may be omitted.
\end{prop}
{\em Proof:}
The relations that are added by the Seifert-van-Kampen theorem are 
$v_{1} u_{2} = 1$, $v_{2} u_{3} = 1$, \dots, $v_{n-1} u_{n} = 1$, $v_{n} u_{1} = 1$. Now as $v_{n} u_{n} v_{n-1} u_{n-1}\dots v_{1}u_{1} = 1$ it is clear that any of the $n$ relations can be omitted. 
\qed 

\section{The representation space} \label{representation space}
The space $R(K;\bf{i})$ consists of representations $\rho \in \text{Hom}(G(K),SU(2))$ such that $\rho(m) \sim \i$, where $m$ is a preferred element (or a set of preferred elements in the case of a link, and the condition of being conjugated to $\i$ is satisfied for each element of the set). Any element in $SU(2)$ with zero trace has order 4 and square $-1$. Therefore, any representation $\rho \in R(K,\bf{i})$ factors through the group 
\begin{equation}
	G(K)_{m,\bf{i}} := (G(K) \times \Z/2) / \langle m^2 (-1) \rangle \ ,
\end{equation}
where $(-1)$ denotes a generator of the group $\Z/2$, and where $\langle m^2 ( -1) \rangle$ denotes the normal subgroup generated by $m^2 (-1)$. Notice that we have a group epimorphism $G(K) \to G(K)_{m,\i}$ induced by the map that sends an element $z \in G(K)$ to the image of $(z,1) \in G(K) \times \Z/2$ in the quotient $(G(K) \times \Z/2) / \langle m^2 (-1) \rangle$. In fact, if we denote the equivalence classes by brackets, then we have 
$[(m^{-2},1)] = [(1,-1)]$, and so the image of the element $(1,-1) \in G(K) \times \Z/2$ lies in the image of the above map $G(K) \to G(K)_{m,\i}$, and therefore the image of any generator of $G(K) \times \Z/2$ is in the image of this map. 

This observation simplifies the description of $R(K,\bf{i})$ considerably. We denote by $\overline{R}(G(K)_{m,\i})$ the representation space $\Hom(G(K)_{m,\i},SU(2))$. Suppose a representation $\rho \in \Hom(G(K)_{m,\i},SU(2))$ satisfies $\rho((-1)) = \rho(m^2) = +1$. This implies that $\rho(m) = \pm 1$. As the knot group $G(K)$ is normally generated by a meridian $m$, this implies that $\rho$  is one of the two possible central representations. These are clearly both isolated in $\overline{R}(G(K)_{m,\i})$. We define
$R(G(K)_{m,\i})$ to be the subspace of non-central representations. 

\begin{prop}\label{alternative description}
The canonical homomorphism $\pi: G(K) \to G(K)_{m,\i}$ induces a homeomorphism
\begin{equation*}
	\pi^* : R(G(K)_{m,\i}) \to R(K;\i) \ ,
\end{equation*}
with both representation spaces seen as subspaces of $\Map(G(K)_{m,\i},SU(2))$, respectively $\Map(G(K),SU(2))$, and with these mapping spaces topologised by the compact-open topology determined by the standard topology on $SU(2)$ and the discrete topology on the groups.   
\end{prop}
{\em Proof: } Any element in $R(K;\i)$ is in the image of $\pi^*$ by construction of $G(K)_{m,\i}$, and as $\pi$ is surjective the map $\pi^*$ is injective likewise. It is an easy matter to check continuity and openness of the map $\pi^*$ directly from the definition of the compact-open topology. 
\qed

\begin{remark}
	The group $G(K)_{m,\i}$ is closely related to the orbifold fundamental group $O(K)= G(K)/\langle m^2 \rangle$ of the 2-fold branched cover $M_2(K)$ of the knot complement branched along the knot, considered as a $\Z/2$-orbifold. In fact, the group $G(K)_{m,\i}$ is a central extension of $O(K)$ by $\Z/2$. On the other hand, the orbifold fundamental group $O(K)$ contains the fundamental group of $M_2(K)$ as a subgroup of index $2$, see for instance \cite[Section 10.6]{Kawauchi}. This also explains the observations made in Section \ref{parallels}. 
\end{remark}

After these general considerations we shall return out attention to pretzel knot or links. 

\begin{prop} For the pretzel knot (or link) $K=P(p_{1}, \dots, p_{n})$, and for $m$ denoting the meridian (respectively a set consisting of one meridian per link component), we have a presentation of $G(K)_{m,\bf{i}}$ given by 
\begin{equation*}
\begin{split}
	\langle s_{1}, \dots, s_{n}| & s_1^4=1, s_1^2 = s_2^2 = \dots = s_n^2, \\
	&
	(s_{1}s_{2})^{p_{1}} = (s_{2} s_{3})^{p_{2}} = \dots = (s_{n-1} s_{n})^{p_{n-1}} = (s_{n}s_{1})^{p_{n}} \rangle .
\end{split}
\end{equation*}
\end{prop}
{\em Proof:} In the presentation of Proposition \ref{presentation knot group} any generator is a meridian $m$. We claim that the group $G(K)_{m,\i}$ has a presentation given by
\begin{equation} \label{presentation 2}
\begin{split}
\langle s_{1}, \dots, s_{n}, -1 | & [-1,s_{i}] = 1, (s_{i})^{2} = -1, i = 1, \dots, n, (-1)^{2} = 1, \\
	&
	(s_{1}s_{2})^{p_{1}} = (s_{2} s_{3})^{p_{2}} = \dots = (s_{n-1} s_{n})^{p_{n-1}} = (s_{n}s_{1})^{p_{n}} \rangle \ .
	\end{split}
\end{equation}
In fact, the relations in the first line follow directly from the definition of $G(K)_{m,\i}$ as a quotient of $G(K) \times \Z/2$. 
The remaining relations of the presentation \ref{presentation knot group} simplify because any of the elements $s_{i}$ now satisfies $s_{i}^{-1} = - s_{i}$. Therefore, in the above notation, the meridians $u$ and $v$ expressed in terms of $s$ and $t$ simplify to the expressions
\begin{align*}
	u & = (st)^p s^{-1} \ \text{ \ and  \ } \\
	 v & = (st)^{p} t^{-1} \text{ \ \ \  for $p$ odd, and  } \\
	u &= (st)^p \  s \text{\ \ \  and \ } \\
	 v & = (st)^{p} \ t \text{ \ \ \ for $p$ even \ ,} 
\end{align*}
no matter of what the sign of $p$ is. The relation $v_1 u_2 = 1$ then becomes $(s_1 s_2)^{p_1} = (s_2 s_3)^{p_2}$, independently of the sign and parity of $p_1$ and $p_2$, and likewise for the remaining ones. Furthermore, the generator $(-1)$ in the presentation (\ref{presentation 2}) may be omitted, yielding the claimed presentation. 

\qed

The elements $\rho \in R(P(p_{1}, \dots, p_{n});\bf{i})$ fall into two classes, depending on whether $\rho((s_{1}s_{2})^{p_{1}}) = \pm 1$ or not, or in other words, according to whether $\rho((s_{1}s_{2})^{p_{1}})$ is central in $SU(2)$ or not. As we shall see, if this is not the case, then the representation $\rho$ is {\em binary dihedral}, which by definition means that it factors through a subgroup of $SU(2)$ that is conjugated to 
\begin{equation*}
	Pin(2) = S^{1} \coprod  \j \cdot S^{1} \ ,
\end{equation*}
where $SU(2)$ is seen as the unit quaternions, $S^{1} \subseteq \C = \langle 1,\i \, \rangle \subseteq \H$ the unit complex numbers, and $\j \cdot S^{1} \subseteq \langle \, \j, \k \, \rangle \subseteq \H$ the circle of unit complex numbers multiplied by $\bf{j}$, lying entirely in the space spanned by $\bf{j}$ and $\bf{k}$. 

Before we proceed, we shall note a useful formula: Let ${\bf v} = v_{1} \, \i + v_{2} \, \j + v_{3} \, \k$ and ${\bf w} = w_{1} \, \i + w_{2} \, \j + w_{3} \, \k$ be purely imaginary quaternions. Then we have
\begin{equation}\label{multiplication of imaginary quaternions}
	{\bf v} \cdot {\bf w}  = - \langle {\bf v},{\bf w} \rangle + {\bf v} \times {\bf w} \ ,
\end{equation}
where 
\[
{\bf v} \times {\bf w} = (v_{2} w_{3} - v_{3} w_{2}) \, \i + (v_{3} w_{1} - v_{1} w_{3} ) \, \j 
	+ (v_{1} w_{2} - v_{2} w_{1}) \, \k \ , 
\]
and where $\langle -, - \rangle$ denotes the standard scalar product.
As the notation suggests, this corresponds to the usual `cross-product' in $\R^{3}$. In particular, if ${\bf v}$ and ${\bf w}$ are linearly independent the vector ${\bf v} \times {\bf w}$ is perpendicular to the plane spanned by ${\bf v}$ and ${\bf w}$. 

\subsection{A conjugacy class fixing condition}
The following Lemma is useful for fixing representants of conjugacy classes of representations. Recall also that the surjective map $S(\text{Im}(\H)) \times [0,\pi] \to SU(2)$, given by $({\bf z}, \alpha) \mapsto e^{{\bf z} \alpha}$ via the exponential map is injective when restricted to $S(\text{Im}(\H)) \times (0,\pi)$, and maps $S \times \{ \pi \}$ to $-1$ and $S \times \{ 0 \}$ to $+1$. 
 
\begin{lemma}\label{conjugacy condition} 
Let $G$ be a group, and let $x$ and $y$ be elements of $G$. Let $\rho: G \to SU(2)$ be a representation such that  $\rho(x)$ is not in the centraliser of $\rho(y)$ (in particular, $\rho$ is non-abelian). Then there are precisely two representations $\rho'$, both conjugated to $\rho$, such that 
\begin{align*}
\rho'(x)  = e^{\j \alpha} \text{ \ \ and \ \ } \rho'(y) = e^{{\bf z} \beta} \ 
\end{align*}
with $\alpha, \beta \in (0,\pi)$, and such that ${\bf z} \in S({\text Im}( \H))$ lies in the plane $\langle \, \j, \k \, \rangle$. The two are distinguished according to the sign of the inner product $\langle \, \k , {\bf z} \rangle \neq 0$ with $\j$. These two are related by conjugation with $\j$, or in other words, related by a rotation of $\z$ around the axis $\j$ with angle $\pi$. 
\end{lemma}
We notice for our further applications that if $\rho(x)$ and $\rho(y)$ are required to have trace zero in $SU(2)$ then we have $\rho'(x) = \j$ and $\rho'(y) = \bf \z$ may also be written as $\rho'(y) = \j e^{\i \alpha'}$ with $\alpha' \in (0,\pi)$ if we require the inner product condition, or with $\alpha' \in (0,2 \pi) \setminus \{ \pi \}$ without the inner product condition. \\

{\em Proof:} 
The action of $SU(2)$ on itself by conjugation consists precisely in the adjoint action on its Lie algebra when seen through the exponential map, and as such it factors through $SO(\su(2))=SO(3)$. 
Up to conjugation we may assume that $\rho(x) = e^{\bf{j} \alpha}$ because the action of $SU(2)$ on the 2-sphere $S({\text Im}( \H)) $ is clearly transitive. This assumption does not yet fix $\rho$ up to conjugation. In fact, we have $c \, \j \, c^{-1} = \j$ if and only if $c \in SU(2)$ is of the form 
$c= w + y \, \j$ with $w,y \in \R$ and $w^{2} + y^{2} = 1$. The set of these elements is a 1-dimensional circle, and conjugating with an element $c$ of this circle yields a rotation around the axis $\j$ in $\langle \i, \j, \k \rangle \cong \R^3$. Therefore, we may assume that $\rho(y)=e^{{\bf z} \beta}$ is such that ${\bf z}$ lies in the $\langle \j, \k \rangle$--plane. Our assumption implies ${\bf z} \neq \j$, and so there are two possible choices of rotations around $\j$ with this requirement on ${\bf z}$, according to whether $\langle \, \k , \bf z\, \rangle <0$ or $\langle \, \k , \bf z\, \rangle > 0$. \\
\hspace{1cm} \qed

\subsection{The case $\rho((s_{1}s_{2})^{p_{1}}) = \rho((s_{2} s_{3})^{p_2}) = \dots = \rho((s_n s_1)^{p_n}) = +1$} This class may have binary dihedral representations and may and usually does contain representations that are not binary dihedral. 

\begin{lemma}\label{perpendicularity}
	Let $\x, \y, \z \in S^2\subseteq \langle{\i,\j,\k \rangle}$, and suppose $\y = \x  \cdot e^{\z \alpha}$ for some $\alpha \in [0,2 \pi]$. Suppose $\x \neq \pm \y$, or equivalently $\alpha \not \equiv 0 \modulo{\pi}$, or equivalently $\sin(\alpha) \neq 0$. Then $\z$ must be perpendicular to $\x$. 
\end{lemma}
{\em Proof:}
As usually one proves that $e^{\z \alpha} = \cos(\alpha) + \z \sin(\alpha)$. Therefore we see that 
\begin{equation}
\x \cdot e^{\z \alpha} = \x \, \cos(\alpha) + \x \cdot \z \, \sin(\alpha) \ .
\end{equation}
As we were assuming $\x \neq \pm \y$ the angle $\alpha$ cannot be $0$ or $\pi$, and so $\sin(\alpha) \neq 0$. 
By the above formula (\ref{multiplication of imaginary quaternions}) we see that this is purely imaginary if and only if $\x$ and $\z$ are perpendicular.
\qed

We are now able to describe the space of all representations up to conjugacy of the pretzel knots or links $P(p_1, \dots, p_n)$. A connection component is reminiscent of the possible configuration space of a mechanical linkage on $S^2 = S({\text Im}( \H))$: It's the space of ordered subsets of points $(\z_1, \dots , \z_n)$ on $S^2$ with (usually) the first two of them fixed by the conjugacy fixing condition, and with the distances of two consecutive points fixed, but with the liberty of `moving around' otherwise. For $n \geq 4$ this yields positive dimensional families in this representation space up to conjugacy.

\begin{prop}\label{+1 case} The set of conjugacy classes of representations 
$\rho \in R(K;\i)$ such that $\rho((s_{1}s_{2})^{p_{1}}) = \rho((s_{2} s_{3})^{p_2}) = \dots = \rho((s_n s_1)^{p_n}) = +1$ is bijective to the ordered subsets $(\z_1, \dots, \z_n)$ of points $\z_1, \dots, \z_n \in S({\text Im}( \H))$, $\z_1 = \rho(s_1), \dots, \z_n = \rho(s_n)$, such that 
\begin{enumerate}
\item the distance between $\z_i$ and $\z_{i+1}$ is given by $\alpha_{i,i+1}  \in [0,\pi]$ satisfying the congruence
\begin{equation*}
p_i \, \alpha_{i,i+1} \equiv p_i \, \pi  \modulo{2\pi} \ ,
\end{equation*}
for $i = 1, \dots, n$, with $n+1=1$ understood, and
\item these points satisfy the following `conjugacy class fixing condition':
\begin{equation*}
	\z_1 =  \j \, , \ \ \ \z_{j+1} = \z_{j} \cdot e^{\i \, \alpha_{j,j+1}}, 
\end{equation*}
for $j=1, \dots, l$, where $l$ is the smallest integer such that $\alpha_{l,l+1} \not \in  \{ 0, \pi \} $, if it exists, or $l=n$ if not. 
\end{enumerate}
Furthermore, a representation $\rho$ determined by an n-tuple $(\alpha_{12}, \alpha_{23}, \dots,
\alpha_{n1})$ is non-abelian unless all angles $\alpha_{i,i+1}$ are equal to $0$ or $\pi$. 
It is binary dihedral if and only if all points $\z_i$ lie on the great circle lying in the $\langle \,  \j, \k \,  \rangle$--plane. Reflection on this plane induces an involution on the space of conjugacy classes of representations with fixed point set precisely the binary dihedral representations. In particular, conjugacy classes of non-abelian representations that are not binary dihedral come in pairs.  
\end{prop}

{\em Proof:} 
The conjugacy fixing condition follows immediately from Lemma \ref{conjugacy condition} above applied to $x = s_1$ and $y = s_{l+1}$ if $\rho$ is non abelian. 
Also by that Lemma we may 
assume that $\rho(s_1) = \j$, and $\rho(s_2) = \j \cdot e^{\i \alpha_{12}}$, with angle $\alpha_{12} \in [0,\pi]$.  Notice that $\alpha_{12}$ is the distance between $\rho(s_1)$ and $\rho(s_2)$ on the 2--sphere $S({\text Im}( \H))$ with its standard metric. Because  $\rho((s_1 s_2))^{p_1} = (-1)^{p_1} e^{\i \alpha_{12} p_1}$ must be equal to $(+1)$, we have the condition that
\begin{equation}\label{condition +1}
	p_1 \, \alpha_{12} \equiv p_1 \pi \modulo{2 \pi} \ .
\end{equation}
We assume for simplicity that $\alpha_{12} \not \in \{ 0, \pi \}$, so that the conjugacy class is already fixed. 
Next we may write $\rho(s_3) = \rho(s_2) \cdot e^{\z_{23} \, \alpha_{23}}$, where $\z_{23} \in S({\text Im}( \H))$ is a purely imaginary quaternion of unit norm. It must be perpendicular to $\z_2 = \rho(s_2) \in S({\text Im}( \H))$ if $\alpha_{23} \not \in \{ 0, \pi \}$ by the above Lemma \ref{perpendicularity}. Similarly to above, the angle $\alpha_{23} \in [0,\pi]$ must satisfy the congruence 
\begin{equation*}
	p_2 \, \alpha_{23} \equiv p_2 \pi  \modulo{2 \pi} \ .
\end{equation*}
We notice that for given angle $\alpha_{23}$ there is a circle of possibilities for the choice of $\rho(s_3)$, parametrised by the circle of elements $\z_{23}$ which are pependicular to $\z_2=\rho(s_2)$, as long as $\alpha_{23}$ is different from $0$ and $\pi$.

This process continues inductively, and the last congruence to satisfy is
\begin{equation*}
	p_n \, \alpha_{n1} \equiv p_n \pi \modulo{2 \pi} \ , 
\end{equation*}
with $\alpha_{n1} \in [0,\pi]$ and $\z_{n1} \in \S $ now such that $\rho(s_1)= \rho(s_n) \cdot e^{\z_{n1} \, \alpha_{n1}}$. In particular, having $\rho(s_n)$ fixed there is only one  possibility of choosing $\z_{n1} \in S({\text Im}( \H))$ -- instead of a whole circle -- as we have to `come back' to $\rho(s_1)$ that was already fixed. 
\qed
\begin{remark}
As the author learned from the referee's report, the observation that there is an involution on the representation space modulo conjugation, with the binary dihedral representations as fixed points, is a general phenomenon. In fact, the group $H^1(G;\Z/2) \cong \Hom(G;\Z/2)$ can be interpreted as the central representations $G \to SU(2)$, and this space acts on the representation space $\Hom(G,SU(2))$ and so on its quotient by conjugation. This action is so that a central representation $\epsilon: G \to \Z/2$ maps a representation $\rho \in \Hom(G,SU(2))$ to the product $\epsilon \cdot \rho$. With a little effort one can check that a representation $\rho$ is a fixed point of this action if and only if $\rho$ is binary dihedral. 
\end{remark}

\subsection{The case $\rho((s_{1}s_{2})^{p_{1}}) = \rho((s_{2} s_{3})^{p_2}) = \dots = \rho((s_n s_1)^{p_n}) = -1$}
This case is entirely analogous to the preceeding one, and the corresponding statement is given by 
\begin{prop}\label{-1 case} The set of conjugacy classes of representations 
$\rho \in R(K;\i)$ such that $\rho((s_{1}s_{2})^{p_{1}}) = \rho((s_{2} s_{3})^{p_2}) = \dots = \rho((s_n s_1)^{p_n}) = -1$ is bijective to the ordered subsets $(\z_1, \dots, \z_n)$ of points $\z_1, \dots, \z_n \in S({\text Im}( \H))$ such that 
\begin{enumerate}
\item the distance between $\z_i$ and $\z_{i+1}$ is given by $\alpha_{i,i+1}  \in [0,\pi]$ satisfying the congruence
\begin{equation*}
p_i \, \alpha_{i,i+1} \equiv (p_i + 1) \pi  \modulo{2 \pi } \ ,
\end{equation*}
for $i = 1, \dots, n$ with $n+1=1$ understood, and
\item these points satisfy the following `conjugacy class fixing condition':
\begin{equation*}
	\z_1 =  \j \, , \ \ \ \z_{j+1} = \z_{j} \cdot e^{\i \, \alpha_{j,j+1}}, 
\end{equation*}
for $j=1, \dots, l$ where $l$ is the smallest integer such that $\alpha_{l,l+1} \neq 0$, if it exists, or $l=n$ if not. 
\end{enumerate}
Furthermore, a representation $\rho$ determined by an n-tuple $(\alpha_{12}, \alpha_{23}, \dots,
\alpha_{n1})$ is non-abelian unless all angles $\alpha_{i,i+1}$ are equal to $0$ (the case $\alpha_{i,i+1} =
\pi$ is excluded by the congruences to satisfy). 
It is binary dihedral if and only if all points $\z_i$ lie on the great circle lying in the $\langle \, \j, \k \, \rangle$--plane. Reflection on this plane induces an involution on the space of conjugacy classes of representations with fixed point set precisely the binary dihedral representations. In particular, conjugacy classes of non-abelian representations that are not binary dihedral come in pairs.  
\end{prop}
{\em Proof:} The difference to the preceeding situation is that now we must have  $-1 = \rho((s_1 s_2))^{p_1} = (-1)^{p_1} e^{\i \alpha_{12} p_1}$ so instead of the condition (\ref{condition +1}) above, we now must have
\begin{equation}\label{condition -1}
	p_1 \, \alpha_{12} \equiv (p_1+1) \pi \modulo{2 \pi} \ , 
\end{equation} 
and so on. 
\qed

\subsection{The case $\rho((s_{1}s_{2})^{p_{1}}) = \rho((s_{2} s_{3})^{p_2}) = \dots = \rho((s_n s_1)^{p_n}) \neq \pm 1$}\label{beta not zero}
As we will see, all representations in this class are binary dihedral: The fact that the distinguished element $(s_1 s_2)^{p_1}= \dots = (s_n s_1)^{p_n}$ is not central in $SU(2)$ forces the images $\rho(s_1), \dots, \rho(s_n)$ of the meridional generators to lie on a great circle.  

\begin{prop}\label{beta case} The set of conjugacy classes of representations 
$\rho \in R(K;\i)$ such that $\rho((s_{1}s_{2})^{p_{1}}) = \rho((s_{2} s_{3})^{p_2}) = \dots = \rho((s_n s_1)^{p_n})$ is conjugate to $e^{\i \beta}$ with $\beta \notin \{0, \pi\}$ is in one-to-two correspondance with the ordered subsets $(\z_1, \dots, \z_n)$ of points $\z_1, \dots, \z_n \in S({\text Im}( \H))$ such that 
\begin{equation*}
	\z_1 = \j \, , \ \ \ \z_{i+1} = \z_{i} \cdot e^{\i \, \alpha_{i,i+1}}, 
\end{equation*}
with the angle $\alpha_{i,i+1}  \in [0,2 \pi]$ satisfying the congruence
\begin{equation*}
p_i \, \alpha_{i,i+1} \equiv \beta \modulo{2 \pi} \ ,
\end{equation*}
for $i = 1, \dots, n$ with $n+1=1$ understood. 

All these representations are binary dihedral and non-abelian.
\end{prop}
{\em Proof:}
Under the assumption we have $\rho(s_{i}) \neq \pm \rho(s_{i+1})$ modulo $n$. Up to conjugation we may assume $\rho(s_{1}) = \bf{j}$ and $\rho(s_{2}) = \j\, e^{\i \alpha}$, with $\alpha \neq 0 \ (\text{mod} \, \pi)$. 

We therefore have $\rho(s_{1} s_{2}) = (-1) \, e^{\i \alpha}$. The image of this element under $SU(2) \to SO(3)$ is given by rotation by $2 \alpha$ around the $\i$ axis, where we consider $\R^{3}$ as the span of $\i,\j,\k$ inside $\H$. By assumption, $e^{\i \beta} := \rho((s_{1} s_{2})^{p_{1}}) \neq \pm 1$, so $\rho((s_{2} s_{3}))^{p_{2}}$ must be a non-trivial rotation around the same axis, the one spanned by $\i$. By the formula (\ref{multiplication of imaginary quaternions}) we therefore see that $\rho(s_{3})$ must lie in the plane $\langle \, \j, \k \, \rangle$ as well, and so may be written as $\rho(s_{3}) = \rho(s_{2}) \, e^{\i \alpha'}$. Inductively, we see that all elements $\rho(s_{i})$ are of the form $\j \, e^{\i \alpha}$ for some angle $\alpha \in [0,2 \pi]$.

\qed
By the methods of Section \ref{arithmetics} below we can see explicitely that there are only finitely many possible values for $\beta$ in the formula $e^{\i \beta} = \rho((s_{1} s_{2})^{p_{1}}) = \dots = \rho((s_n s_1)^{p_n})$. However, this also follows from  \cite[Theorem 10]{Klassen}.

\begin{prop}
  If the number of strands of the pretzel knot is $n=3$ then the representation space modulo conjugation $\mathscr{R}(K;\i):= R(K;\i)/SU(2)$ only consists of isolated points. 

  If the number of strands of the pretzel knot is $n \geq 4$, and if $\mathscr{R}(K;\i)$ contains a representation that is not binary dihedral, then it contains connection components of strictly positive dimension. 
\end{prop}
{\em Sketch of proof:} 
For the claim regarding $n = 3$ strands notice that a non binary dihedral representation is determined by a triangle on the 2--sphere with the lengths of its sides and two edge points fixed and with the lengths satisfying the arithmetic conditions given above. There are only finitely many such representations for each knot. There are only finitely many binary dihedral representations by Klassen's result. 

For $n \geq 4$ the conjugacy classes in $\mathscr{R}(K;\i)$ will in general come in positive dimensional families -- because one may `move' the n-gon even if two consecutive points of it are fixed. 
\qed 
\begin{remark}
	We will compute the Zariski tangent spaces to $\mathscr{R}(K;\i)$ for a large class of pretzel knots or links $K$ below, and our results will be coherent with the preceeding observation. We expect that `generically' $\mathscr{R}(P(p_1,\dots,p_n;\i)$ contains components of dimension $n-3$ for $n \geq 3$. For $n=4$ the 1-dimensional components are given by 1-spheres, and for higher $n$ the author expects the $n-3$ dimensional components to be spheres of that dimension. 
\end{remark}

\subsection{Orbits of the conjugacy action}
A representation $\rho \in \Hom(G,SU(2))$ is called {\em irreducible} if there is no proper subspace $V$ of
$\C^2$ that is invariant under $\rho$, in the sense that $\rho(g) V = V$ for all $g \in G$. Otherwise it is
called reducible. It is not hard to see that a representation into $SU(2)$ is irreducible if and only if
it is non-abelian. 

Let us consider the action of $SU(2)$ on $\Hom(G,SU(2))$ given by conjugation. A representation $\rho$ is
irreducible if and only if its stabiliser $\Gamma_\rho \subseteq SU(2)$ is equal to the centre $\Z/2$. The
trivial representation and representations with image inside the centre of $SU(2)$  have stabiliser $SU(2)$,
and a reducible representation that acts non-trivially on a proper subspace has stabiliser isomorphic to
$U(1)$. Therefore, the orbit $[\rho]$ of an irreducible representation is isomorphic to $SU(2)/\Z/2 \cong
\mathbb{RP}^3$, and the orbit of a reducible representation that acts non-trivially on a proper subspace is
homeomorphic to $SU(2)/U(1) \cong S^2$. In the situation of the representation spaces $R(K;\i)$ that we
consider the reducible representations with stabiliser $SU(2)$ do not appear. 

\subsection{Abelian representations}
The reducible/abelian representations inside the space $R(K;\i)$ are quite easily described for a general knot. 
\begin{prop}\label{abelians}
	Let $K=$ be a knot or link. Then there are, up to conjugation, precisely $2^{\abs{K}}$ abelian representations in $R(K;\i)$ with $\abs{K}$ denoting the number of components of $K$. 
\end{prop}
{\em Proof:}
Any abelian representation $\rho: \pi_1(S^3 \setminus K) \to SU(2)$ factors through the abelianisation $H^1(S^3 \setminus K;\Z) \cong \Z^{\abs{K}}$, where $\abs{K}$ denotes the number of components of $K$, and each meridian represents a generator. For a knot there is therefore just one representation up to conjugacy such that the meridian is mapped to an element conjugate to $\i$. In the case of a link the requirement that $\rho(m_1) = \i$, where $m_1$ denotes a meridian to the first component of $K$, implies that $\rho(m_i) = \pm \i$ for meridians to the other components $m_2, \dots, m_{\abs{K}}$. The claim follows. 
\qed
It is easy to see that in the case of a pretzel knot or link $P(p_1, \dots, p_n)$ we have the following: If all $p_i$ are odd then our result is a knot in the case that $n$ is odd and is a 2-component link in the case that $n$ is even. Otherwise the number of components is equal to the number of $p_i$'s that are even.

\section{Arithmetic properties}\label{arithmetics}
We show that under a simple arithmetic condition on the numbers $p_1, \dots, p_n$  irreducible representations with certain `degenerate properties' may be avoided, with the best possible situation for $n=3$.  
\begin{prop}\label{arithmetic 1}
	Suppose the numbers $p_1, \dots, p_n$ are pairwise coprime. Then there is no representation $\rho$ determined by angles 
	\begin{equation*}
	(\alpha_{1,2}, \dots, \alpha_{n,1})
	\end{equation*}
with only one or two of the angles $\alpha_{i,i+1}$ not in the set $\{0,\pi\}$. In particular, if $n=3$, then any $\rho$ that is non-abelian must have angles $\alpha_{12}, \alpha_{23}, \alpha_{31}$ none of which lies in $\{0,\pi\}$.
\end{prop}
{\em Proof:}
We only show that there is no representation with only two of the angles not in the set $\{0,\pi\}$, leaving the easier other case as an exercise to the reader. 

	Notice that if $\rho$ is not such that $\rho((s_{1} s_{2})^{p_{1}}) = \dots = \rho((s_n s_1)^{p_n}) = \pm 1$ then none of the angles $\alpha_{i,i+1}$ may be in the set $\{0,\pi\}$. So we may assume $\rho$ to be such that this equation holds, and therefore in the situation of Proposition \ref{+1 case} or \ref{-1 case} above.
	
Suppose $\alpha_{i,i+1}, \alpha_{j,j+1}$ with $i \neq j$ are such that they both do not lie in $\{0, \pi \}$, and all others do. These two angles must satisfy
\begin{align*}
	& p_i \ \frac{\alpha_{i,i+1}}{\pi} - p_i \equiv p_j \frac{\alpha_{j,j+1}}{\pi} - p_j \equiv 0 \modulo{2\Z} \ \ \mbox{or} \\
& p_i \ \frac{\alpha_{i,i+1}}{\pi} - p_i \equiv p_j \frac{\alpha_{j,j+1}}{\pi} - p_j \equiv 1 \modulo{2\Z} \ , 
\end{align*}
depending on whether we have $\rho((s_1s_2)^{p_1}) = \pm 1$. 
Therefore we may assume that $\alpha_{i,i+1}/\pi = k_i/p_i$, and $\alpha_{j,j+1}/\pi = k_j / p_j$ with numbers $k_i, k_j \in \Z$ (the parity of these numbers depends on the case in which we are, and on the parity of $p_i$ and $p_j$, but this is not important in the following argument). 
On the other hand it is easy to see that the points $\z_i = \rho(s_i)$ map to precisely two points under $S({\text Im}( \H)) = S^2 \to \mathbb{RP}^2$, and in particular lie on a geodesic circle joining two non-antipodal points on $S^2$, of length $2 \pi$. Therefore we must have the congruence
\begin{equation*}
\frac{\alpha_{i,i+1}}{\pi} + \frac{\alpha_{j,j+1}}{\pi} \equiv 0 \modulo{\Z} . 
\end{equation*} 
This implies that there is an integer $n \in \Z$ such that
\begin{equation*}
	k_i p_j + k_j p_i = n \, p_i p_j \ , 
\end{equation*}
and so $p_j$ divides $k_j p_i$. As $p_i$ and $p_j$ are assumed to be coprime, $p_j$ must divide $k_j$. But this contradicts that $\alpha_j = \pi \frac{k_j}{p_j}$ does not lie in $\{0, \pi\}$. 
\qed 

Our next arithmetic result concerns binary dihedral representations. 
\begin{prop}\label{arithmetic 2}
	Suppose the  numbers $p_1, \dots, p_n$ are pairwise coprime, and suppose $\rho$ is a binary dihedral representation. Then either $\rho$ is abelian or we must have (up to conjugation)
\begin{equation*}
\rho((s_{1} s_{2})^{p_{1}}) = \dots = \rho((s_n s_1)^{p_n}) = e^{\i \, \beta}
\end{equation*}
with $\beta \notin \{ 0, \pi \}$. 
In other words, the situation in Section \ref{beta not zero} above is the one that must occur for all non-abelian binary dihedral representations under this assumption.
\end{prop}
{\em Proof:}
Suppose $\rho$ is binary dihedral and $\frac{\beta}{\pi}$ is an integer. We may suppose that $\rho(s_1) = \j$ and 
\begin{equation*}
\rho(s_{i+1}) = \rho(s_i) \, e^{\i \, \alpha_{i,i+1}} \ ,
\end{equation*}
with angle $\alpha_{i,i+1} \in [0,2\pi]$, for $i=1, \dots, n$ and $n+1 = 1$ understood. This implies that the sum of the angles must be a multiple of $2 \pi$, 
\begin{equation}\label{condition binary dihedral}
	\frac{\alpha_{1,2}}{\pi} + \dots + \frac{\alpha_{n-1,n}}{\pi} + \frac{\alpha_{n,1}}{\pi} \equiv 0 \modulo{2 \Z} .
\end{equation}
In addition, we must have the congruences 
\begin{equation*}
p_i \frac{\alpha_{i,i+1}}{\pi} \equiv 0  \modulo{\Z}
\end{equation*}
for $i=1, \dots, n$. Putting $\alpha_{i,i+1} /\pi = k_i / p_i$ with $k_i \in \Z$, $i=1, \dots, n$, inserting this in equation (\ref{condition binary dihedral}), and multiplying this equation by $p_1 \cdot \dots \cdot p_n$ we see that $p_j$ divides $k_j \, p_1 \cdot \dots \cdot \widehat{p}_{j} \cdot \dots \cdot p_n$, with the hat on $\widehat{p}_j$ indicating that this factor is omitted. By the condition on pairwise coprimeness we see that $p_j$ must in fact divide $k_j$, and this for $j = 1, \dots, n$. As a consequence, each angle $\alpha_{i,i+1}$ must be $0$ or $\pi$, and so the representation $\rho$ is abelian. 
\qed 

\section{Non-degeneracy conditions}\label{non deg}
The local structure of the representation variety  $\Hom(G,SU(2))$ of a discrete group $G$ was first studied by Weil, see \cite{Weil1, Weil}. A presentation $\langle g_1, \dots, g_n | r_1, \dots, r_m \rangle$ of $G$, identifies the space $\Hom(G,SU(2))$
 with $F^{-1}(1,\dots,1)$, where $F: SU(2)^n \to SU(2)^m$ is given by the polynomials determined by the relations that the generators have to satisfy.  

%

\begin{definition}\label{Zariski tangent space} Let $\rho \in \Hom(G,SU(2))$. A map $\xi: G \to \su(2)$ is called a cocycle at $\rho$ if one has
\begin{equation}\label{cocycle condition}
\xi(gh) = \xi(g) + \Ad_{\rho(g)} \xi(h) 
\end{equation}
for any $g,h \in G$. 
An element $\zeta \in \su(2)$ defines a coboundary $\zeta^\#: G \to \su(2)$ at a representation $\rho$ by the formula
\begin{equation}\label{coboundary condition}
\zeta^\#(g) = \zeta - \Ad_{\rho(g)} \zeta 
\end{equation}
for $g \in G$. 
\end{definition}
To motivate this definition, let $\rho$ be a representation. A deformation into nearby representations
$\rho_t$ may be written as
$\rho_t(g) = \rho(g) + t \, \xi(g) \rho(g) + o(t)$. It is then easy to check that the
requirement of $\rho_t$ to be a group homomorphism implies for $\xi$ to be a cocycle. 

Coboundaries are cocycles. They are inifinitessimal deformations of $\rho$ that are induced by conjugating $\rho$ by elements of $SU(2)$. In fact, the  
coboundary $\zetadiez$ associated to $\zeta \in \su(2)$ is obtained as the derivative of $t \mapsto e^{t\zeta} \rho(g) e^{-t \zeta}$ at $0$. 

The space of cocycles at $\rho$ is denoted by $Z^1(G;\su(2)_\rho)$, the space of coboundaries by $B^1(G;\su(2)_{\rho})$, and the quotient is denoted by $H^1(G;\su(2)_\rho)$. As the notation suggests, this is isomorphic to  the first cohomology group of $K(G;1)$ with twisted
coefficients in $\su(2)$ defined by the adjoint action of $G$ on $\su(2)$ determined by $\rho$, see \cite{Weil}. Weil proves that if the space $H^1(G;\su(2)_\rho) $ vanishes, then the representation $\rho$ has a neighbourhood in $\Hom(G,SU(2))$ (with the compact open topology) all of which elements are conjugated to $\rho$. In other words, the conjugacy class of $\rho$ is isolated in $\Hom(G,SU(2))/SU(2)$. This suggests to call a representation {\em non-degenerate} if $H^1(G;\su(2)_\rho) = 0$, and this is the definition we will take for our purpose. 

This definition is not always suitable: For instance, if $G$ is a knot group, then one has $\dim H^1(G;\su(2)_\rho) \geq 1$ for cohomological reasons, and it can be shown that in the case that this dimension is precisely 1, the space $\Hom(G,SU(2))/SU(2)$ has the structure of a 1-dimensional smooth  manifold in a neighbourhood of an irreducible representation $\rho$, see for instance \cite[Lemma 2 and Proposition 1]{HeusenerKlassen}. 

In general, the space $H^1(G;\su(2)_\rho)$ is referred to as the Zariski tangent space of $\Hom(G,SU(2))/SU(2)$ at $\rho$. However, it can appear that a representation $\rho$ is isolated in $\Hom(G,SU(2))/SU(2)$ but that its Zariski tangent space is non-trivial. This is well-known, but in the cases we consider examples will even appear below.

In this section we shall study the local structure of $R(K;\i)$, seen as $R(G(K)_{m,\i}) = \Hom(G(K)_{m,\i},
SU(2))$ by Proposition \ref{alternative description} above. In particular, a non-degenerate representation
$\rho \in R(K;\i)$ may well have deformations not coming from the action by conjugation when seen as an
element of the bigger representation space $\Hom(G(K),SU(2))$ (without the assumption that meridians are
mapped onto trace-free matrices). 

\begin{lemma}\label{cocycles at -1}
If $\rho$ is a non-abelian representation then a coboundary $\xi$ at $\rho$ necessarily satisfies $\xi((-1)) = 0$. 
\end{lemma}
{\em Proof:} In fact, by the cocycle condition we must have
\begin{equation*}
	\xi(g\, (-1)) = \xi(g) + \Ad_{\rho(g)} \xi((-1)) = \xi((-1)) + \Ad_{\rho((-1))} \xi(g) = \xi((-1) \, g)  
\end{equation*}
for all $g \in G(K)_{m,\i}$ as $(-1)$ commutes with all elements in $G(K)_{m,\i}$. Clearly the endomorphism $\Ad_{\rho((-1))}$ is the identity. Therefore $\Ad_{\rho(g)} \xi((-1)) = \xi((-1))$ for all $g \in G(K)_{m,\i}$. If $\rho$ is non-abelian this implies $\xi((-1))=0$. \qed

\begin{lemma}\label{cocycles at meridians}
If $\xi$ is a cocycle at $\rho$ then for any element $h$ that is (conjugated to) a meridional element $m$ the element $\xi(h) \in \su(2)$ must be perpendicular to $\rho(h) \in S^2 \subseteq \langle \, \i, \j, \k \, \rangle = \su(2) = \R^3$. 
\end{lemma}
{\em Proof:} Indeed, as we have the relation $h^2=-1$ in $G(K)_{m,\i}$ the cocycle condition together with the preceeding Observation implies 
\begin{equation*}
0 = \xi(h^2) = (1 + \Ad_{\rho(h)}) \, \xi(h) \ .
\end{equation*}
By assumption $\Ad_{\rho(h)}$ is a rotation by $\pi$ around the axis $\rho(h)$. Consequently $\xi(h)$ must lie in the plane annihilated by $(1+\Ad_{\rho(h)})$ which is precisely the plane of elements perpendicular to $\rho(h)$. \qed

\begin{prop}\label{non-deg dihedral}
Let $K=P(p_1, \dots, p_n)$ be a pretzel knot or link with non-zero determinant $\abs{\Delta_K(-1)}$ (without any arithmetic assumption on the $p_i$'s). 
Let $\rho \in R(K;\i)$ be a non-abelian representation which is binary dihedral, and which satisfies $\rho((s_1 s_2)^{p_1}) = \dots = \rho((s_n s_1)^{p_n}) \neq \pm 1$. Then the cohomology group $H^1(G_{K;\i},\su(2)_\rho)$ vanishes.
\end{prop}
{\em Proof:}
Up to conjugation we may assume that $\rho(s_1) = \j$ and $\rho(s_{i+1}) = \rho(s_{i}) \, e^{\i \, \alpha_{i,i+1}}$ with angles $\alpha_{i,i+1} \in [0,2 \pi]$ for $i=1, \dots, n$ and $n+1 = 1$ understood. The assumption implies that there is an angle $\beta \in [0,2 \pi]$, different from $0$ and $\pi$, such that $p_i \,  \alpha_{i,i+1} \equiv \beta \modulo{2 \pi \Z}$ and $\rho((s_1 s_2)^{p_1}) = \dots = \rho((s_n s_1)^{p_n}) = (-1) \,  e^{\i \beta} $ . 
Suppose $\xi$ is a cocycle at $\rho$. Notice that we have
\begin{equation}\label{B}
	\xi((s_i \, s_{i+1})^{p_i}) = \underbrace{(1 + \Ad_{\rho(s_i s_{i+1})} + \dots 
						+ \Ad_{\rho(s_i s_{i+1})}^{p_i -1})}_{=: B_{i,i+1}} \, \xi(s_i \, s_{i+1}) \ , 
\end{equation}
with $\Ad$ denoting the adjoint action of $SU(2)$ on the Lie algebra $\su(2)$. 
\begin{lemma}\label{automorphism}
	The endomorphism $B_{i,i+1}$ of $\su(2)$ is an automorphism. 
\end{lemma} 
{\em Proof of the Lemma:}
Notice that we have the equation
\begin{equation*}
	1 - \Ad_{e^{\i  \beta}} = B_{i,i+1} (1 - \Ad_{\rho(s_i s_{i+1})} ) 
\end{equation*}
for any $i= 1, \dots, n$. Now $\Ad_{e^{\i \beta}}$ is rotation by angle $2 \beta \notin \{ 0, 2 \pi \}$ around the $\i$--axis. Therefore $(1- \Ad_{e^{\i \beta}})$ has the subspace spanned by $\i \in \su(2) \cong \R^3$ as its kernel and maps the whole space onto the plane perpendicular to $\i$. Therefore $B_{i,i+1}$ must have rank at least $2$. Likewise, $\Ad_{\rho(s_i s_{i+1})}$ is a non-trivial rotation around the $\i$--axis, and so $(1-  \Ad_{\rho(s_i s_{i+1})})$ maps the plane perpendicular to $\i$ onto itself. Therefore $B_{i,i+1}$ must map the plane perpendicular to $\i$ onto itself. On the other hand, $B_{i,i+1}$ restricted to the subspace spanned by $\i$ is just multiplication by the number $p_i$, and so $B_{i,i+1}$ is an automorphism.
\qed 

For an element $\lambda \in \su(2) = \langle \, \i, \j, \k \, \rangle $ we denote by  $\lambda \|_{\i}$  its projection onto the subspace generated by $\i$. We will prove the Proposition by showing that the space of cocycles $Z^1_\rho(G_{K,\i};\su(2))$ is equal to the space of coboundaries $Z^1_\rho(G_{K,\i};\su(2))$. 
\\

Let us first observe that if $\xi$ and $\xi'$ are two coboundaries at $\rho$ such that $\xi(s_1) = \xi'_(s_1)$ and $\xi(s_2) = \xi'(s_2)$ then in fact $\xi = \xi'$. To see this, note that $(s_1 s_2)^{p_1} = (s_2 s_3)^{p_2}$ and the assumption that $\xi$ and $\xi'$ coincide on $s_1$ and $s_2$  implies that
\[
\xi((s_2 s_3)^{p_2}) = \xi'((s_2 s_3)^{p_2}) \ ,
\]
or equivalently, by the above formula (\ref{B}) and the Lemma \ref{automorphism},
\[
\xi(s_2 s_3) = \xi'(s_2 s_3) \ .
\]
By the cocycle formula this of course implies that $\Ad_{\rho(s_2)} \xi(s_3) = \Ad_{\rho(s_2)} \xi'(s_3)$, and so $\xi(s_3) = \xi'(s_3)$. Inductively, we show in the same way that $\xi(s_4) = \xi'(s_4), \ \dots \ ,$ $  \xi(s_n) = \xi'(s_n)$. 

Therefore it is enough to show that for any cocycle $\xi$ there is an element $\zeta$ of the Lie algebra $\su(2)$, and so a coboundary $\zeta^\#$, such that $\xi(s_1) = \zetadiez(s_1)$ and $\xi(s_2) = \zetadiez(s_2)$. To show this, we characterise the possible values of $(\zetadiez(s_1),\zetadiez(s_2))$ more precisely. Consider the map
\begin{equation*}
\begin{split}
	\su(2) & \to \rho(s_1)^{\perp} \oplus \rho(s_2)^{\perp}
 \\
	\zeta & \mapsto ((1 - \Ad_{\rho(s_1)}) (\zeta) , (1 - \Ad_{\rho(s_2)}) (\zeta)) 
		= (\zetadiez(s_1),\zetadiez(s_2)) \ . 
\end{split}
\end{equation*}
We will show now that this map has rank 3 and an element $(\lambda,\mu) \in \rho(s_1)^{\perp} \oplus \rho(s_2)^{\perp}$ is in the image of this map if and only if 
\begin{equation*}
	\langle \lambda, \i \rangle = \langle \mu , \i \rangle \ , 
\end{equation*}
or equivalently, if and only if $\lambda$ and $\mu$ have identical $\i$--component. 
 For this notice first that $\z_1=\rho(s_1)$ and $\z_2=\rho(s_2)$ have been assumed to lie perpendicular to $\i$. We decompose $\su(2)$ as $\langle \,  \i \, \rangle \oplus \langle \, \i \, \rangle ^{\perp}$, writing an element $\zeta = \zeta_\| + \zeta_\perp$ correspondingly. We find that
\begin{equation*}
\begin{split}
	\langle \, (1- \Ad_{\z_1})(\zeta) - (1-\Ad_{\z_2})(\zeta) , \i \, \rangle 
	& = \langle \, \Ad_{\z_2} \zeta + \Ad_{\z_1} \zeta , \i \, \rangle \\
	& = \langle \, \Ad_{\z_2} \zeta_\perp + \Ad_{\z_1} \zeta_\perp , \i \, \rangle \\
	& = 0 \ , 
\end{split}
\end{equation*}
because $\Ad_{\z_1}$ and $\Ad_{\z_2}$ certainly preserve the subspace $\langle \, \i \, \rangle ^{\perp}$. This implies that we do indeed have $\langle \, \lambda - \mu , \i \, \rangle = 0$ for any element $(\lambda, \mu)$ in the image of the above map. It is injective because the assumption implies that $\rho(s_1) \neq \pm \rho(s_2)$. As the image space $\rho(s_1)^\perp \oplus \rho(s_2)^\perp$ has rank 4, we  therefore see that an element $(\lambda,\mu)$ of this space that satisfies $\langle \, \lambda - \mu , \i \, \rangle = 0$ must lie in the image. 
\\

We finish the proof by showing that a cocycle $\xi$ at $\rho$ satisfies $\langle \, \xi(s_1) - \xi(s_2) , \i \, \rangle = 0$ and therefore is a coboundary by the above. 
Recall that the automorphisms $B_{i,i+1}$ introduced above respects the splitting $\langle \,  \i \, \rangle \oplus \langle \, \i \, \rangle ^{\perp}$ of $\su(2)$, and that $B_{i,i+1}$ restricted to the span of $\i$ is just multiplication by $p_i$. Now expressing the cocyle $\xi$ at the element $(s_1 s_2)^{p_1} = \dots = (s_n s_1)^{p_n}$ in the $\i$--direction and using this fact we just obtain
\begin{equation*}
p_1 (\xi(s_1)_{\|_{\i}} - \xi(s_2)_{\|_{\i}}) = p_2 (\xi(s_2)_{\|_{\i}} - \xi(s_3)_{\|_{\i}}) = \dots = p_n (\xi(s_n)_{\|_{\i}} - \xi(s_1)_{\|_{\i}}) \ . 
\end{equation*}
Let $c \in \R$ be this number. We only have to show that $c = 0$. 
However, notice that the sum
\[
	\sum_{i=1}^{n} (\xi(s_i)_{\|_{\i}} - \xi(s_{i+1})_{\|_{\i}})
\]
is just $0$, as this is a cyclic sum because $s_{n+1} = s_1$ was understood. Multiplying this sum with the product $p_1 \dots p_n$, we therefore find
\begin{equation*}
\begin{split}
	0 & = p_1 \cdot \dots \cdot p_n \left( \sum_{i=1}^{n} (\xi(s_i)_{\|_{\i}} - \xi(s_{i+1})_{\|_{\i}}) \right) \\
	& = ( \sum_{i=1}^{n} p_1 \cdots \hat{p_i} \cdots p_n ) \cdot c \\
	& = \pm \Delta_{K}(-1) \cdot  c \ ,
\end{split}
\end{equation*}
where $\Delta_{K}(-1)$ is just the determinant of our Pretzel knot or link $K$. By our assumption this number is non-zero, and therefore we may conclude that $c=0$, and so that $\xi$ is a coboundary. 
\qed

\begin{prop}
Let $K=P(p,q,r)$ be a pretzel knot and let $\rho \in R(K;\i)$ be a non-abelian representation with $\rho((s_1 s_2)^{p_1}) = \dots = \rho((s_n s_1)^{p_n}) = \pm 1$, and
 which is not binary dihedral. Then the cohomology group $H^1_\rho(G_{K;\i},\su(2))$ vanishes: The Zariski tangent space of $R(K;\i)$ at $\rho$ is equal to the space of coboundaries determined by $\rho$. 
\end{prop}
{\em Proof:} 
Up to conjugation we may assume that $\rho(s_1) = \j$ and $\rho(s_{2}) = \j \, e^{\i \, \alpha_{12}}$ with angle $\alpha_{12} \in [0,2 \pi]$. Suppose $\xi$ is a cocycle at $\rho$. The conclusions of Observation \ref{cocycles at -1} and \ref{cocycles at meridians} remain valid in this situation. Therefore, the dimension of the space of cocylces at $\rho$ is at most $6$. However, instead of Lemma \ref{automorphism} we now have:
\begin{lemma} \label{B in the case +1 or -1}
Suppose that we have $\rho((s_i s_{i+1})^{p_i}) = \pm 1$ and that $\rho(s_i) \neq \pm \rho(s_{i+1})$. Then the endomorphism $B_{i,i+1}$ of $\su(2)$ defined in equation (\ref{B}) above has rank $1$. More precisely, we have $B_{i,i+1} = p_i \, \Pi_{\rho(s_i) \times \rho(s_{i+1})}$, where $\Pi_{\rho(s_i) \times \rho(s_{i+1})}$ is the projection onto the space spanned by $\rho(s_i) \times \rho(s_{i+1}) \in \su(2)$. 
\end{lemma}
{\em Proof:}
As we now have $\Ad_{\rho((s_i s_{i+1})^{p_i}} = \id_{\su(2)} $ we can conclude that
\begin{equation*}
 0 = B_{i,i+1} (1 - \Ad_{\rho((s_i s_{i+1}))}) \ . 
\end{equation*}
As by assumption  $\rho(s_i) \neq \pm \rho(s_{i+1})$  we know that $\Ad_{\rho((s_i s_{i+1}))})$ is a non-trivial element in $SO(\su(2))$. Therefore the image of $(1- \Ad_{\rho((s_i s_{i+1}))})$ must have rank $2$, and so the kernel of $B_{i,i+1}$ must at least contain the 2-dimensional subspace of $\su(2)$ that is perpendicular to the rotation axis $ \langle \, \rho(s_i) \times \rho(s_{i+1}) \, \rangle $ of $\Ad_{\rho(s_i s_{i+1})}$. On the other hand, it is immediate from the definition of $B_{i,i+1}$ that it is given by multiplication by $p_i$ when restricted to the 1-dimensional subspace $ \langle \, \rho(s_i) \times \rho(s_{i+1}) \, \rangle $.
\qed

Because of $(s_1 s_2)^p = (s_2 s_3)^q = (s_3 s_1)^r$ the cocycle must satisfy
\begin{equation*}
\begin{split}
  B_{12} \, \xi(s_1 s_2) - B_{23} \, \xi(s_2 s_3) & = 0 \\
  B_{23} \, \xi(s_2 s_3) - B_{31} \, \xi(s_3 s_1) & = 0 \ .
\end{split}
\end{equation*}
By assumption the representation $\rho$ is not binary dihedral, and so the three axes
$\rho(s_1) \times \rho(s_{2}), \, \rho(s_2) \times \rho(s_{3})$ and $ \, \rho(s_3) \times \rho(s_{1})$ are pairwise linearly independent. Therefore, the last equations are equivalent to 
\begin{equation*}
\begin{split}
  B_{12} \, \xi(s_1 s_2) & = 0 \\
  B_{23} \, \xi(s_2 s_3) & = 0 \\
  B_{31} \, \xi(s_3 s_1) & = 0 \ .
\end{split}
\end{equation*}
Equivalently, the element $(\xi(s_1),\xi(s_2),\xi(s_3))$ lies in the kernel of the linear map 
$L: \su(2)^3 \to \su(2)^3$ given by 
\begin{equation*}
 (\xi_1,\xi_2,\xi_3) \mapsto \underbrace{\begin{pmatrix} B_{12} & 0 & 0 \\ 0 & B_{23} & 0 \\ 0 & 0 & B_{31} \end{pmatrix}}_{:=B} \underbrace{ \begin{pmatrix} 1 & \Ad_{\rho(s_1)} & 0 \\
				0 & 1 & \Ad_{\rho(s_2)} \\
				\Ad_{\rho(s_3)} & 0 & 1 
		\end{pmatrix}
		}_{:=C}
		\begin{pmatrix} \xi_1 \\ \xi_2 \\ \xi_3
		\end{pmatrix} \ .
\end{equation*}

The linear map defined by $B$ has rank $3$ by the preceeding Lemma.
On the other hand the map $C$ is an automorphism. In fact, elementary line operations show that the determinant of $C$ is just equal to the determinant of the endomorphism \[
1 + \Ad_{\rho(s_3)} \circ \Ad_{\rho(s_1)} \circ \Ad_{\rho(s_3)} = 1 + \Ad_{\rho(s_3) \rho(s_1) \rho(s_2)} 
\]
 of $\su(2)$. For this it is sufficient to show that $-1$ cannot be an eigenvalue of $\Ad_{\rho(s_3) \rho(s_1) \rho(s_2)} \in SO(\su(2))$, so cannot be a rotation by angle $\pi$, and this itself is equivalent to showing that $\rho (s_3) \rho(s_1) \rho(s_2)$ can not be a purely imaginary quaternion. By the above assumption $\rho(s_1) \rho(s_2) = - e^{\i \alpha_{12}}$. Let $\rho(s_3)=: \z \in S(\text{Im}(\H))$. We now compute, using the above formula (\ref{multiplication of imaginary quaternions}), 
 \begin{equation*}
 	\rho(s_3)  \rho(s_1) \rho(s_2) = \z \cdot (- e^{\i \alpha_{12}}) 
	= - \langle \, \z , \i \, \rangle \, \sin(\alpha_{12}) 
	+ \z \times \i \, \sin(\alpha_{12}) + \z \, \cos(\alpha_{12}) \ . 
 \end{equation*}
This is purely imaginary if and only if $\langle \, \z , \i \, \rangle = 0$, which is true if and only if $\rho$ is binary dihedral. As this was excluded by assumption, we conclude that $-1$ is not an eigenvalue of $\Ad_{\rho(s_3) \rho(s_1) \rho(s_2)}$, and so $C$ indeed is an automorphism. 
Consequently, the homomorphism $L$ has rank $3$. 

Let $H_i \subseteq \su(2)$ be the subspace in which the cocycle $\xi(s_i)$ must lie according to Lemma \ref{cocycles at meridians}. 
Now let's consider the restriction of of $L$ to to the subspace $H_1 \oplus H_2 \oplus H_3$. We claim that this restriction still has rank 3. 

Indeed, let $\chi_{12} \in \text{im}(B_{12})$ be an element in the 1-dimensional image of $B_{12}$. We claim that there is an element $\zeta_1 \in H_1$ such that
\begin{equation} \label{to solve}
	(\chi_{12},0,0) = L\, (\zeta_1,0,0) \ . 
\end{equation}
To see this, notice that for any $\zeta_1 \in H_1$ one has 
\begin{equation*}
	L\, (\zeta_1,0,0) = (B_{12} \,  \zeta_1 , 0 , B_{31} \, \Ad_{\rho(s_3)} \, \zeta_1) \ . 
\end{equation*}
The endomorphism $\Ad_{\rho(s_3)}$ preserves the kernel of $B_{31}$. Thus  (\ref{to solve}) has solutions for non-trivial $\chi_{12}$ if there is an element $\zeta_1 \in H_1$ which is in the kernel of $B_{31}$, and which projects non-trivially onto $\text{im}(B_{12})$, or equivalently, that does not lie in the kernel of $B_{12}$ either. But the kernels of $B_{12}$ and $B_{31}$ only coincide if $\rho(s_1), \rho(s_2)$ and $\rho(s_3)$ all lie in the same plane in $\text{Im}(\H)$. This can only occur if $\rho$ is binary dihedral which is excluded by assumption. Finally we have to convince ourselves, that we can pick an element in the kernel of $B_{31}$ which is not in the kernel of $B_{12}$ and which also lies in $H_1$. But this follows because the intersection of the kernel of $B_{12}$ and the kernel of $B_{31}$ is precisely the span of $\rho(s_1)$. But $H_1 = \langle \rho(s_1) \rangle^{\perp}$ is precisely the orthogonal complement. So $H_1$ and $\ker(B_{31})$ intersect in a 1-dimensional subspace that does not lie in the kernel of $B_{12}$. 

Likewise one shows that for an arbitrary $\chi_{23} \in \text{im}(B_{23})$ and for an arbitrary $\chi_{31} \in \text{im}(B_{31})$ the element $(0,\chi_{23},0)$ respectively $(0,0,\chi_{31})$ is in the image of $L$ restricted to $H_1 \oplus H_2 \oplus H_3$. Therefore, this restriction $L|_{H_1 \oplus H_2 \oplus H_3}$ has rank 3, and its kernel is also of rank 3. 

%
%

So the fact that $(\xi(s_1),\xi(s_2),\xi(s_3))$ must lie in the kernel of $L$ then implies that the space of cocycles at $\rho$ is 3-dimensional, and so is equal to the space of coboundaries. 
\qed 

\begin{remark}
	The method of the above proof applies to the situation of pretzel knots or links with $n \geq 4$ strands as well. Generically, the equivalent of the map $L$ above, restricted to $H_1 \oplus \dots \oplus H_n$, will have rank $n$. Therefore, we expect the dimension of the Zariski tangent space $H^1(P(p_1, \dots, p_n);\su(2)_\rho)$ to be equal to $n-3$ for a generic non binary dihedral representation $\rho$. 
	
\end{remark} 

\begin{theorem}
	Let $K=P(p,q,r)$ be a pretzel knot with $p,q,r$ pairwise coprime. Then the twisted cohomology group $H^1(G(K)_{m, \i};\su(2)_\rho)$ vanishes at any representation $\rho \in R(K;\i)$. 	 
\end{theorem}
{\em Proof:} 
    As we are in the case of a knot, the claim is true at the (conjugacy class) of the unique abelian representation in $R(K;\i)$. At any representation $\rho$ which is not binary dihedral the preceeding Proposition applies. At any non-abelian binary dihedral representation $\rho$ the arithmetic assumption implies that $\rho((s_1 s_2)^{p_1}) \neq \pm 1$ by Proposition \ref{arithmetic 2} above. 
Therefore, Proposition \ref{non-deg dihedral} comes to bare with the desired conclusion. 
\qed

\begin{remark}
	If the pretzel knot $P(p,q,r)$ admits a binary dihedral representation $\rho$ with $\rho((s_1 s_2)^{p_1}) = \dots = \rho((s_3 s_1)^{p_3}) = \pm 1$ then this representation is degenerate in the sense that the associated cohomology group $H^1(G(P(p,q,r))_{m, \i};\su(2)_\rho)$ is non-vanishing. For instance, the pretzel knot $P(3,3,3)$ is such a knot. Nonetheless, these representations are isolated in the representation space modulo conjugacy $\mathscr{R}(P(p,q,r);\i)$. 
\end{remark} 

{\em Sketch of proof:} In this case the maps $B_{12}, B_{23}$ and $B_{31}$ with the terminology from above all have equal 2-dimensional kernels $M$. Again, with the above terminology, the space $M \cap H_i$ is 1-dimensional for $i=1, \dots, 3$. 
 Any element $\xi$ with $\xi(s_i) \in M \cap H_i$ chosen arbitrarily for each $i$ defines a cocycle at $\rho$, so the space of all $\xi$ with this property is 3-dimensional. However, only a 1-dimensional subspace of cocycles with this property is given by coboundaries. 

The claim about isolatedness of the points in $\mathscr{R}(P(p,q,r);\i)$ follows from the results in Section \ref{representation space} above. 
\qed 

\begin{remark}
It is interesting to notice that instead of computing the Zariski tangent space explicitely one may also get non-degeneracy results by studying the topology of the 3-manifold one gets from taking the double branched cover of $S^3$, branched along the knot, see the corresponding results of \cite{HeusenerKlassen}.
\end{remark}

\section{The pretzel knots $P(p,q,r)$}\label{P(p,q,r)}
We get a complete picture of the representation spaces $R(K;\i)$ of the pretzel knots or links $K=P(p,q,r)$. First of all, notice a few pathologies: A pretzel knot or link $P(0,q,r)$ is a sum of of a $(2,q)$ torus knot or link and a $(2,r)$ torus knot or link. If one of $\abs{p}, \abs{q}, \abs{r}$ is $1$ then we are left with a knot or link with bridge number at most $2$, as is easy to see. A 2-bridge knot or link only has binary dihedral representations in $R(K;\i)$ .
We do not have in mind these pathological knots or links. 
For simplicity, the following Proposition is only stated for knots, but a corresponding statement for links is self-suggesting. 
\begin{prop}\label{summary for representation space}
Let $K$ be the pretzel knot $P(p,q,r)$ with all of  $p, q, r$ 
different from $0$ or $\pm 1$, and such that these numbers are pairwise coprime . Then the 
representation space $R(K;\i)$ is isomorphic to the disjoint union
\begin{equation*}
 S^2 \coprod \left( \coprod_{I} \mathbb{RP}^3 \right) \ , 
\end{equation*}
where the finite index set $I$ parametrises the conjugacy classes of all non-abelian 
representations. Among these there are $(\abs{\Delta_{P(p,q,r)}(-1)}-1)/2$ many binary dihedral ones, as well as $I-(\abs{\Delta_{P(p,q,r)}(-1)}-1)/2 > 0$ many
non-binary dihedral representations that are described in Proposition \ref{+1 case} and
\ref{-1 case} above. 
\end{prop}
{\em Proof:}
The fact that there is a single orbit homeomorphic to $S^2$ follows from Lemma \ref{abelians} above. That there are only finitely many isolated orbits homeomorphic to $\mathbb{RP}^3$ follows from the results of the preceeding sections. The number of non-abelian binary dihedral conjugacy classes was determined by Klassen \cite{Klassen}. So we only have to show that there are non-binary dihedral representations. 

We may without loss of generality assume that $\abs{r} \geq \abs{q} \geq \abs{p}$. We show that there is a non-abelian representation $\rho$ as in Proposition \ref{+1 case}. To do so, we will fix the first two angles $\alpha_{pq}$ and $\alpha_{qr}$ satisfying the congruence $p \, \alpha_{pq} \equiv p \pi \modulo{2 \pi}$ and $q \, \alpha_{qr} \equiv q \pi \modulo{2 \pi}$.  It will remain to show that we can find a distance $\alpha_{rp} \in (0,\pi)$ which satisfies the congruence $r \, \alpha_{rp} \equiv r \pi \modulo{2 \pi}$, and such that the triangle inequality 
\begin{equation}\label{triangle inequality}
	\abs{\alpha_{pq}-\alpha_{qr}} \leq \alpha_{rp} \leq \alpha_{pq} + \alpha_{qr}
\end{equation}
holds. 

There are three cases to consider. First, we assume $p$ and $q$ are both odd. 
We just pick $\alpha_{pq} = \pi/\abs{p}$ and $\alpha_{qr} = \pi/\abs{q}$, both satisfying the required congruences for $\alpha_{pq}$ and $\alpha_{qr}$ in this case. 
The triangle inequality (\ref{triangle inequality}) is then equivalent to  
\begin{equation*}
	\frac{\abs{q}-\abs{p}}{\abs{pq}} \leq \frac{\alpha_{rp}}{\pi} \leq \frac{\abs{q} + \abs{p}}{\abs{pq}} \ .
\end{equation*}
The interval $[ \frac{\abs{q}-\abs{p}}{\abs{pq}}, \frac{\abs{q} + \abs{p}}{\abs{pq}}]$ has length $\frac{2}{\abs{q}}$. But as we have assumed $\abs{r} \geq \abs{q}$, there certainly are at least two integer multiples of $\frac{\pi}{\abs{r}}$ inside this interval, and the choice of one of them for the number $\alpha_{rp}/\pi$ will satisfy the congruence $r \, \alpha_{rp} \equiv r \, \pi  \modulo{2 \pi}$. 

Next, assume that $p$ is even and $q$ is odd. We choose $\alpha_{pq} = 2\pi/\abs{p}$ and $\alpha_{qr} = \pi/\abs{q}$.
The triangle inequality (\ref{triangle inequality}) is then equivalent to  
\begin{equation*}
	\frac{2\abs{q}-\abs{p}}{\abs{pq}} \leq \frac{\alpha_{rp}}{\pi} \leq \frac{2\abs{q} + \abs{p}}{\abs{pq}} \ .
\end{equation*}
Again, the interval $[ \frac{2\abs{q}-\abs{p}}{\abs{pq}}, \frac{2\abs{q} + \abs{p}}{\abs{pq}}]$ has length $\frac{2}{\abs{q}}$, and again one may pick a multiple of $\frac{\pi}{\abs{r}}$ inside this interval for the number $\alpha_{rp}/\pi$, so that it satisfies the required congruence. 

If finally $q$ is even and $p$ is odd, we choose $\alpha_{pq} = \pi/\abs{p}$ and $\alpha_{qr} = 2\pi/\abs{q}$. The triangle inequality (\ref{triangle inequality}) is then equivalent to  
\begin{equation*}
	\frac{\abs{\; \abs{q}-2\abs{p}\; }}{\abs{pq}} \leq \frac{\alpha_{rp}}{\pi} \leq \frac{{\abs{q} + 2\abs{p}}}{\abs{pq}} \ .
\end{equation*}
The interval $[\frac{\abs{\; \abs{q}-2\abs{p}\; }}{\abs{pq}}, \frac{\abs{q}+2\abs{p}}{\abs{pq}}]$ now has length $\frac{4}{\abs{q}}$ or $\frac{2}{\abs{p}}$. In both cases, one may pick a multiple of $\frac{\pi}{\abs{r}}$ as $\alpha_{rp}/\pi$ inside this interval, and this satisfies the required congruence. 

In each situation, it follows from Proposition \ref{arithmetic 1} above that the corresponding representation is not binary dihedral.
\qed 
It is interesting that we get the following Corollary. The result is not new, see \cite{BoileauZieschang}. 
\begin{corollary}\label{bridge number}
	Let $K$ be the pretzel knot $P(p,q,r)$ with the numbers $p,q$ and $r$ pairwise coprime and non-zero. Then $P(p,q,r)$ has bridge number $3$ if and only if all of $\abs{p},\abs{q},\abs{r}$ are strictly bigger than $1$. 
\end{corollary}
{\em Proof:} Inspection of the diagram shows that the bridge number can at most be $3$. If any of $\abs{p},\abs{q},\abs{r}$ is equal to $1$ then it is easy to see that the knot is 2--bridge. If non of these number is $1$, then the preceeding Proposition implies the existence of non binary dihedral representations. However, a 2--bridge knot $K$ only has representations in $R(K;\i)$ that are binary dihedral. \qed

\begin{remark}
	It seems that the Proposition \ref{summary for representation space} and its Corollary are true without the arithmetic assumption of the numbers $p,q$ and $r$ being pairwise coprime. However, without this assumption, it looks like the proof would require a little more work in order to sort out solutions where the above inequalities (\ref{triangle inequality}) are both strict. 
\end{remark}

As an example, we shall now compute the non-abelian representations of
any of the knots $P(3,5,7), P(-3,5,7), P(-3,-5,7), P(3,-5,7), \dots$ (all combinations of signs may occur) that are not binary dihedral. As a matter of notation, we shall write $\overline{\alpha}_{i,i+1}:= \alpha_{i,i+1} / \pi$ for the angles occuring in Proposition \ref{+1 case} and \ref{-1 case} above. These have to satisfy the congruences in these Propositions, and as distances between the points $\z_1=\rho(s_1), \z_2=\rho(s_2)$, and $\z_3=\rho(s_3)$, they have to satisfy the triangle inequality. 
\\

The representations $\rho$ with $\rho((s_{1}s_{2})^{3}) = \rho((s_{2} s_{3})^{5} = \dots = \rho((s_3 s_1)^{7}) = +1$ are determined by the following table, that first lists all possible combinations of angles satisfying the congruencies, and then checks the triangle-inequality on each. 
\begin{equation*}
\begin{tabular}{c|c|c|| c|c || c}
$\overline{\alpha}_{12}$ & $\overline{\alpha}_{23}$ & $\overline{\alpha}_{31}$  & $\abs{\overline{\alpha}_{23} - \overline{\alpha}_{31}}$ & $\overline{\alpha}_{23} + \overline{\alpha}_{31}$ & $\Delta \text{--inequality}$
\\
\hline
\hline
1/3 & 1/5 & 1/7 &  2/35 & 12/35 & \text{no} \\
    &     & 3/7 &  8/35 & 22/35 & \text{yes} \\
    &     & 5/7 & 18/35 & 32/35 & \text{no} \\ 
    & 3/5 & 1/7 & 16/35 & 26/35 & \text{no} \\
    &     & 3/7 & 6/35  & 36/35 & \text{yes} \\
    &     & 5/7 & 4/35  & 46/35 & \text{yes.} \\
\end{tabular} 
\end{equation*}
Likewise, the representations $\rho$ with $\rho((s_{1}s_{2})^{3}) = \rho((s_{2} s_{3})^{5} = \dots = \rho((s_3 s_1)^{7}) = -1$ are determined by the following table.

\medskip
\begin{equation*}
\begin{tabular}{c|c|c|| c|c || c}
$\overline{\alpha}_{12}$ & $\overline{\alpha}_{23}$ & $\overline{\alpha}_{31}$  & $\abs{\overline{\alpha}_{23} - \overline{\alpha}_{31}}$ & $\overline{\alpha}_{23} + \overline{\alpha}_{31}$ & $\Delta \text{--inequality}$
\\
\hline
\hline
2/3 & 2/5 & 2/7 &  4/35 & 24/35 & \text{yes} \\
    &     & 4/7 &  6/35 & 34/35 & \text{yes} \\
    &     & 6/7 & 16/35 & 44/35 & \text{yes} \\ 
    & 4/5 & 2/7 & 18/35 & 38/35 & \text{yes} \\
    &     & 4/7 & 8/35  & 48/35 & \text{yes} \\
    &     & 6/7 & 2/35  & 58/35 & \text{yes.} \\
\end{tabular} \ 
\end{equation*}
Each combination of angles that gives rise to a non-abelian representation that is not binary dihedral yields precisely two different conjugacy classes of representations. Therefore, there are in total 18 such conjugacy classes for any of the knots $P(3,5,7), P(-3,5,7), P(-3,-5,7), P(3,-5,7), \dots$. 


\section{Relation to Lin's knot invariant and Heusener and Kroll's generalisation}
Lin has defined a knot invariant, that he denotes $h(K)$, from the representation space $\mathscr{R}(K;\i) = R(K;\i)/SU(2)$ considered here. In the case that all irreducible representations are non-degenerate, and so these are isolated points in $\mathscr{R}(K;\i)$, the number $h(K)$ is just a signed count of these conjugacy classes of irreducible representations. Surprisingly, this is related to the signature of $K$ in the following way \cite{Lin}:
\begin{equation*}
	h(K) = \frac{1}{2} \, \sign(K) \ .
\end{equation*}
The signature of a knot is just the signature of the symmetric bilinear form given by the matrix $V + V^{t}$, with $V$ denoting a Seifert matrix of the knot.

For the pretzel knots $P(p,q,r)$ with $p,q,r$ odd the signature is easily computed. Indeed, a Seifert matrix for these genus-1-knots is given by
\begin{equation*}
	V = \frac{1}{2} \begin{pmatrix} p+q & q+1 \\ q-1 & q+r \end{pmatrix} \ ,
\end{equation*}
see for instance \cite[Example 6.9]{Lickorish}. Depending on the numbers $p,q,r$ the signature is just 2,0 or -2, and so Lin's invariant $h(K)$ equals 1,0 or -1. If all of $p,q,r$ are in addition of the same sign, then the signature is 1 or -1, and in particular Lin's invariant is odd. As the non-abelian non-binary dihedral representations come in pairs, there must be an odd number of binary dihedral representations of which there are $\frac{1}{2} (\abs{\Delta_K(-1)}-1)$ by Klassen's result, thus reflecting the fact that the determinant $\Delta_P(-1)= pq +pr +qr$ is congruent to $3$ modulo $4$ in this case. 

The knot $P(-3,5,7)$ with trivial Alexander polynomial, and therefore no binary dihedral representations at all, has 0 signature. Lin's result then reflects the fact again that the non-binary-dihedral representations come in pairs. 

The pretzel knots with one strand with an even number of crossings have higher genus in general and also higher signature. For instance, the Fintushel-Stern knot $P(-2,3,7)$ has genus 5 and signature 8. It has determinant 1 and so no binary dihedral representations at all. 

In \cite{HeusenerKroll} Heusener and Kroll extend Lin's result to the situation of studying spaces $R^\alpha(K)$ of representations $\rho$ modulo conjugation, such that $\rho(m) \sim e^{\i \alpha} \in SU(2)$. They define an invariant $h^{\alpha}(K)$ as a count of $\mathscr{R}^{\alpha}(K)= R^{\alpha}(K)/SU(2)$, and establish $h^{\alpha}(K) = \frac{1}{2} \sign_K(e^{\i 2 \alpha})$, thereby extending Lin's result. Here $\sign_K: S^1 \setminus \{ 1 \} \to \Z$ is the Levine-Tristram signature function, i.e. $\sign_K(\omega)$ is the signature of the Hermitian form $(1-\omega) V + (1-\overline{\omega}) V^{t}$, where $\omega \in S^1 \setminus \{ 1 \} \subseteq \C $, and $V$ is a Seifert matrix of $K$.

Our notion of non-degeneracy for representation $\rho \in R(K;\i)$ was based on representations of the group $G(K)_{m,\i}$. In particular, it implied that its conjugacy class $[\rho] \in \mathscr{R}(K;\i)$ isolated. Of course, without the restriction that $\rho$ maps the meridian to an element conjugated to $\i$ in $SU(2)$ this doesn't need to be true. In fact, when seen as an element of the representation space $\mathscr{R}(K)$ of all representations of the knot group in $SU(2)$, up to conjugation, this ceases to be true, as follows for instance from Heusener and Kroll's result \cite[Corollary 3.9]{HeusenerKroll}, where it is proved from Lin's result and a continuity property for their extension of Lin's invariant that a knot with non-zero signature must have irreducible representations $\rho$ that have a 1-parameter family of deformations $(\rho_t)$ inside the full space $\mathscr{R}(K)$ coming with a non-trivial family of deformations $\tr(\rho_t(m))$, and so deform out of $\mathscr{R}(K;\i) \subseteq \mathscr{R}(K)$. 

The method of introducing the group $G(K)_{m,\i}$ in order to study more easily the representation space $R(K;\i)$ might be applied to the study of $R^\alpha(K)$ with $e^{\i \alpha}$ another root of unity. However, the simplification is probably largest for the case of a fourth root of unity as considered in this article.

\end{document}